\documentstyle[12pt]{article}

\newtheorem{theorem}{Theorem}[section]

\newtheorem{lemma}{Lemma}[section]
\newtheorem{corollary}{Corollary}[section]
\newtheorem{prop}{Proposition}[section]

\def\k#1{\kern#1em}
\def\Ib#1{{\rm I\kern-.25em#1}}
\def\Ibb#1{{\rm I\kern-.23em#1}}

\def\vcg{\vrule  width.02em height1.4ex depth-.05ex}
\def\C{{\rm\k{.24}\vcg\k{-.26}C}}

\def\R{\Ib R}
\def\N{\Ibb N}

\def\P{\Ib P}

\def \RR {{\rm Res}\,}
\def\({{\rm (}}
\def\){{\rm )}}
\begin{document}
\title{Analytic residues along algebraic cycles
\thanks{AMS classification 
number: 32A27, (32A25, 32C30)}}
\author { by \\ Carlos. A. Berenstein \\
carlos@src.umd.edu\\
 Institute for Systems Research\\
           University of Maryland, MD 20742, USA  \\ 
           Alekos Vidras \\
         msvidras@ucy.ac.cy\\ 
         Department of Mathematics and Statistics\\
           University of Cyprus, Nicosia 1678, Cyprus \\
         Alain Yger \\
         yger@math.u-bordeaux.fr\\
         Department of Mathematics, \\
         University of Bordeaux,  33405 Talence,  France.  }
\maketitle
\newpage
\begin {abstract}
Let $W$ be a $q$-dimensional irreducible algebraic 
subvariety in the affine space ${\bf A}^n_\C$,   
$P_1,...,P_m$ $m$ elements in $\C[X_1,...,X_n]$, and $V(P)$ the set of 
common zeros of the $P_j$'s in $\C^n$. Assuming 
that $|W|$ is not included in $V(P)$, one can attach to 
$P$ a family of non trivial $W$-restricted residual currents 
in $'{\cal D}^{0,k}(\C^n)$, $1\leq k\leq \min(m,n)$, with support 
on $|W|$. These currents 
(constructed following an analytic approach) 
inherit most of the properties that are fulfilled in the case $q=n$.  
When the set $|W| \cap V(P)$ is discrete and $m=q$, 
we prove that for every point $\alpha \in |W|\cap V(P)$
the $W$-restricted analytic residue of a $(q,0)$-form $R d\zeta_I$, $R\in \C[X_1,...,X_n]$,  
at the point $\alpha$ 
is the same as the residue on ${\cal W}$ (completion of $W$ in ${\rm Proj}\,\C[X_0,...,X_n]$)  
at the point $\alpha$ in the sense of Serre 
($q=1$) or Kunz-Lipman ($1<q<n$) of the $q$-differential form $(R/P_1\cdots P_q) d\zeta_I$. 
We will present a restricted version of some affine version of Jacobi's residue formula
and applications of this formula to  
higher dimensional analogues of Reiss (or Wood) relations, 
corresponding to situations where the Zariski closures of $|W|$ and $V(P)$ intersect at infinity
in an arbitrary way. 
\end{abstract}
\newpage
\section{Introduction}
\setcounter{equation}{0}
Let $\Gamma$ be a complete integral curve 
embedded as a closed subscheme in   
${\rm Proj}\, \C[X_0,...,X_n]$ 
and $\C(\Gamma)$ its function field. 
Following the exposition of H\"ubl and Kunz 
of the Serre's approach \cite{huku:gnus}, 
the residue of a meromorphic $(1,0)$-differential 
form $\omega \in \Omega^1_{\C( \Gamma )/\C}$ at the point 
$\alpha \in \Gamma$ is defined as follows~: 
let ${\cal M}_1, \dots ,{\cal M}_d$ be the minimal prime ideals of the completion 
$\widehat {\cal O}_{\Gamma ,\alpha}$ 
of the local ring of $\Gamma$ at $\alpha$ and 
let ${\overline R}_j$, $j=1,...,d$, be the integral closures of the 
"branches" $R_j=\widehat {\cal O}_{\Gamma ,\alpha} /{\cal M}_j$, 
$j=1,\dots, d$, of the curve $\Gamma $ at the point $\alpha$.  
Then ${\overline R}_j$ is isomorphic to some algebra of formal power series 
$\C[[t_j]]$ and in  $\C ((t_j))$ the differential
$(1,0)$-form $\omega$ can be written as
$$
\omega=\sum\limits_{k\geq k_{j}} a_k^j\; t_j^k \,,
$$
where $a_k^j\in \C$, $k\geq k_{j}$, are complex numbers which 
are independent of the parameters $t_j$. Define 
\begin{eqnarray*}
{\rm Res}_{\Gamma,\alpha,{\overline R}_j }\, \omega :=a^{j}_{-1}\,,\quad
{\rm Res}_{\Gamma,{\alpha }}\, \omega:=\sum \limits_{j=1}^d a^{j}_{-1}\, . 
\end{eqnarray*}

\vskip 2mm
\noindent
It was pointed by G. Biernat in \cite{bier1:gnus} that, if $f_1,...,f_n$ 
are $n$ germs of holomorphic functions in $n$ variables 
(with jacobian determinant $J_f\in {\cal O}_n$) such that 
$(f_1,...,f_{n-1})$ define a germ of curve $\gamma$ (with branches parametrized 
respectively by $\varphi_1,...,
\varphi_d$) and ${\rm dim}\, [\gamma\cap \{J_f=0\}]=0$, then, for 
any $h\in {\cal O}_n$, the Grothendieck residue 
$$
{\rm Res}_0\, \Big[{h d\zeta_1\wedge \cdots \wedge d\zeta_n 
\over f_1\cdots f_n}\Big]:= {1\over (2i\pi)^n} \int\limits_{{|f_1|=\epsilon_1}
\atop {\vdots \atop {|f_n|=\epsilon_n}}}
{h d\zeta_1 \wedge \cdots \wedge d\zeta_n \over f_1\cdots f_n} 
$$
(with the orientation for the cycle $\{|f_1|=\epsilon_1\,,\dots\,, |f_n|=\epsilon_n\}$ that 
ensures the positivity of the differential form $d\,  {\rm arg}\, f_1 \wedge \cdots \wedge 
d\, {\rm arg}\, f_n$ on it) 
equals 
$$
\sum\limits_{j=1}^d {\rm Res}_{t=0} \, \Bigg[\Big({f_n' h \over f_n J_f}
\circ \varphi_j\Big) (t)\, dt\Bigg]\,; 
$$
in particular, if $\omega$ denotes the $(1,0)$-meromorphic differential form 
$$
\omega:= {g d\zeta_\alpha \over f_n} \,,\quad g\in {\cal O}_n\,,\ \alpha \in \{1,...,n\}\,,
$$
then
$$
{\rm Res}_0 \, \Big[{ {df_1 \wedge \cdots \wedge df_{n-1}\over f_1\cdots f_{n-1}} \wedge \omega}\Big] 
$$
equals the sum 
\begin{eqnarray}
\sum\limits_{j=1}^d \nu_j {\rm Res}_{\gamma_j,0}\, [\omega]\,,
\end{eqnarray}
where $\gamma_1,...,\gamma_d$ correspond to the irreducible germs of curves attached  
to the isolated primes in the decomposition of $(f_1,...,f_{n-1})$, 
and ${\rm Res}_{\gamma_j,0}\, [\omega]$ is defined 
on the model of the Kunz-H\"ubl residue, this notion being transposed from the 
algebraic context to the analytic one. This suggests a natural relation between 
the approaches developped by Serre-H\"ubl-Kunz and the analytic residue approach 
developped by Coleff-Herrera \cite{coh:gnus} (which precisely allows the transposition of 
the definition of the Grothendieck residue in the complete intersection case to 
the setting of currents). 
\vskip 2mm
\noindent
The analytic approach we use to define restricted residual currents on 
a $q$-dimensional reduced analytic space ${\cal Y}\subset U$, where $U$ is an open 
subset of $\C^n$, will be described in section $2$ as follows~: if $f_1,...,f_m$ are $m$ 
functions holomorphic in $U$, then the map
$$
\lambda \mapsto \Phi_{Y,f}(\lambda):=\|f\|^{2 \lambda} \, [Y]\,,
$$
where $[Y]$ denotes the integration current on $Y=|{\cal Y}|$, can be meromorphically 
continued as a $'{\cal D}^{(n-q,n-q)}(U)$-map. Moreover, for any $k\in \{1,...,m\}$, 
for any ordered subset ${\cal I}\subset \{1,...,m\}$ with cardinal $k\leq \min(q,m)$, the 
analytic continuation of 
$$
\lambda \mapsto \lambda c_k \Phi_{Y,f} (\lambda-k-1) \wedge \overline \partial 
\|f\|^2 \wedge \Big( \sum\limits_{l=1}^k (-1)^{l-1} \overline {f_{i_l}} 
\bigwedge\limits_{{j=1}\atop {j\not=l}}^k \overline {df_{i_j}}\Big)\,, 
$$
where 
$$
c_k:={(-1)^{k(k-1)/2} (k-1)! \over (2i\pi)^k} 
$$
is holomorphic at the origin. Its value at $0$ defines, up to a multiplicative 
constant, a residual regular holonomic $(n-q,n-q+k)$-current which is supported by $Y\cap V(f)$~; 
regular holonomiticity is here understood in the sense of Bj\"ork (\cite{bjo3:gnus}, 
chapter $9$). Properties of such currents are similar 
to those introduced as above in the case $q=n$. Proposition 2.1 will summarize 
the different properties of such restricted residual currents. The main case of interest 
for us will be the case where $m\leq q$ and ${\rm dim} \,(Y\cap V(f)) \leq q-m$, that is 
$f_1,...,f_m$ define a complete intersection in ${\cal Y}$. In this case, the 
restricted residue current corresponding to ${\cal I}:=\{1,...,m\}$ is 
the Coleff-Herrera current on $Y$  
$$
\Big(\bigwedge \limits_{j=1}^m \overline \partial {1\over f_j}\Big) \wedge [Y]
$$          
introduced in \cite{coh:gnus}. It is not surprizing that residual restricted currents 
in such a complete intersection setting obey the transformation law for 
multidimensional residue calculus (\cite{gh:gnus}, chapter 6), which we will prove (and use next) in the 
case $m=q$. If $f_1=P_1,...,f_q=P_q$ are polynomials and $W$ is an affine 
$q$-dimensional algebraic subvariety 
of the affine scheme ${\bf A}^n_\C$ such that ${\rm dim}\, (V(P)\cap |W|)=0$, we will 
prove in section $3$ that the total sum of restricted residues 
$$
\RR 
\left[\matrix{ [W]\wedge Q dX_{i_1}\wedge \cdots \wedge dX_{i_q} \cr
\cr P_1,..., P_q}\right] 
$$
vanishes as soon as the degree of $Q$ is sufficiently small, under a properness 
assumption on the restriction of $(P_1,...,P_q)$ to $|W|$. We will thus transpose to the 
restricted case an Abel-Jacobi formula proved in the case $q=n$ and $W={\bf A}^n_\C$ in 
\cite{vy:gnus}.    
\vskip 2mm
\noindent
Let again $W$ be a $q$-dimensional irreducible algebraic subvariety 
in the affine scheme ${\bf A}^n_\C$   
and $P_1,...,P_q$, $q$ polynomials in $\C[X_1,...,X_n]$ such that 
$|W|\cap V(P)$ is a discrete (hence finite) algebraic set in $\C^n$. Let 
$Q\in \C[X_1,...,X_n]$ and ${\cal I}$ a subset in $\{1,...,n\}$ with 
cardinal $q$. The meromorphic differential form 
$$
\omega:= {Q d\zeta_{i_1}\wedge \cdots \wedge d\zeta_{i_q} \over P_1\cdots P_q} 
$$
induces an element in $\Omega^q_{\C({\cal W})/\C}$, where ${\cal W}$ 
denotes the completion of $W$ in ${\rm Proj}\, \C[X_0,...,X_n]$. 
We will prove in section 4, 
thanks to the algebraic residue theorem in \cite{lip2:gnus} and 
the properties of restricted residual currents that were pointed out in previous 
sections, that 
the residue at a closed point $\alpha$ in 
$|W|\cap V(P)$ (in the sense of H\"ubl or Lipman \cite{lip1:gnus}) of the 
differential form $\omega$ (viewed as an element in $\Omega^q_{\C({\cal W})/\C}$)
equals 
$$
{\rm Res}_{W,\alpha}\, [\omega]:=\Big\langle 
\Big(\bigwedge \limits_{j=1}^q \overline \partial {1\over P_j}\Big) \wedge [W]\,,\, 
\psi Q d\zeta_{i_1}\wedge \cdots \wedge d\zeta_{i_q}\Big\rangle \,,
$$
where $\psi$ denotes a test-function with compact support in some arbitrary small 
neighborhood of $\alpha$, such that $\psi \equiv 1$ near $\alpha$. The result is 
clear when $\alpha$ is a smooth point of $W$, it will follow from the 
algebraic residue formula combined with a perturbation argument in the 
case $\alpha$ is a singular point of $W$. As a consequence of the 
fact that the analytic and algebraic approaches lead to the same restricted residual 
objects, we will extend (also in section $4$) with an algebraic formulation  
to such a restricted context
the affine Jacobi's theorem obtained in the non-restricted case $W={\bf A}^n_\C$ in \cite{vy:gnus}. 
\begin{theorem} 
Let $W$ be a $q$-dimensional irreducible affine algebraic subvariety in ${\bf A}^n_\C$ \($0<q<n$\) 
and $P_1,...,P_q$ be $q$ polynomials in $\C[X_1,...,X_n]$ such that there exists 
strictly positive rational numbers $\delta_1,...,\delta_q$ and two constants $K>0$,
 $\kappa>0$ 
such that~: 
\begin{eqnarray} 
\zeta \in |W|\,,\ \|\zeta\|\geq K \ \Longrightarrow \ \sum\limits_{j=1}^q {|P_j(\zeta)|\over 
\|\zeta\|^{\delta_j}} \geq \kappa\,;
\end{eqnarray} 
then, for any $Q\in \C[X_1,...,X_n]$ such that ${\rm deg}\, Q <\delta_1+\cdots +\delta_q -q$,
for any multi-index $(i_1,...,i_q)$ in $\{1,...,n\}^q$,  
\begin{eqnarray}
\sum\limits_{\alpha\in |W|\cap V(P)} 
{\rm Res}_{W,\alpha}\, \Big[{ Qd\zeta_{i_1}\wedge \cdots \wedge d\zeta_{i_q} 
\over P_1\cdots P_q}\Big] =0\,.
\end{eqnarray} 
\end{theorem} 
We will derive (in sections $4$ and 5) some consequences of this result in the spirit 
of Cayley-Bacharach's theorem or Wood's results \cite{wood1:gnus}. The key point 
here (compare to the framework of \cite{huku:gnus} or \cite{ku2:gnus}) 
is that the properness assumption along $|W|$ (1.2) which is satisfied by the polynomial 
map $P:=(P_1,...,P_q)$ does not imply that the Zariski closures of $|W|$ 
and $V(P_1,...,P_q)$ in $\P^n(\C)$ have an empty 
common intersection on the hyperplane 
at infinity.     

\section{Restricted residual currents}
\setcounter{equation}{0}

We begin this section by recalling some basic facts about currents on analytic  
manifolds, especially integration currents on analytic sets or Coleff-Herrera currents and 
their ``multiplication'' with integration currents~; we will inspire ourselves 
from \cite{bjo:gnus}, \cite{by2:gnus}, \cite{by3:gnus} or \cite{meo:gnus}. 
\vskip 2mm
\noindent
We start with basic facts about integration on a $q$-dimensional irreducible analytic subset 
$Y$ in $U\subset \C^n$ \cite{le:gnus}. The subset $Y_{\rm reg}$ of regular points
of $Y$ is a $q$-dimensional complex manifold. The set of singular points 
$Y_{\rm sing}$ is an analytic subset of $U$ with complex dimension ${\rm dim}\, Y_{\rm sing}<q$. 
Therefore for any smooth $(q,q)$ test form $\phi_{(q,q)} \in {\cal D}^{(q,q)}(U)$, 
one can define the action of the integration current $[Y]$ on $\phi_{(q,q)}$ as  
\begin{eqnarray*}
\langle [Y]\,,\, \phi_{(q,q)}\rangle 
&=&\int _{Y}\phi _{(q,q)}(\zeta,\bar \zeta)= \int _{Y_{\rm reg}}\phi_{(q,q)}(\zeta,\bar \zeta)+
\int_{Y_{\rm sing}}\phi _{(q,q)}(\zeta,\bar \zeta)\\
&=& \int_{Y_{\rm reg}}\phi_{(q,q)}(\zeta,\bar \zeta)\, .  
\end{eqnarray*} 
For ${\rm Re}\, \lambda >0$ and $f_1,...,f_m$ holomorphic in $U$, one can define the 
$(q,q)$-current $\|f\|^{2\lambda} \, [Y]$ by 
$$
\langle \|f\|^{2\lambda}\, [Y]\,,\, \phi_{(q,q)}\rangle \, := 
\int_{Y_{\rm reg}} \|f\|^{2\lambda} \, \phi_{(q,q)}\,.
$$
It is known (\cite{bjo:gnus},\cite{bjo3:gnus}) that this current $[Y]$ is a 
regular holonomic current, which implies, for each point $z_0$ in 
$U\cap Y$, the existence of a Bernstein-Sato relation 
\begin{eqnarray}
{\cal Q}_{z_0} \Big(\lambda, \zeta, \bar \zeta, {\partial \over 
\partial \zeta}, {\partial \over \partial \bar \zeta}\Big)\;  
\Big[ \|f\|^{2(\lambda+1)} \otimes [Y]\Big] =b_{z_0}(\lambda)\, (\|f\|^{2\lambda} \otimes [Y])
\end{eqnarray}
($b_{z_0}\in \C[X]$) 
valid in a neighborhood of $z_0$.  
In fact, this does not follow directly from theorem 3.2.6 in \cite{bjo:gnus} since 
$\|f\|^2$ is a real analytic function (and not a holomorphic one). Nevertheless, the 
existence of Bernstein-Sato relations of the form (2.1) remains valid here since $\|f\|^2$ 
has the particular form 
$$
\|f(\zeta)\|^2=\sum\limits_{j=1}^m f_j(\zeta) \overline {f_j (\zeta)} 
$$
and the integration current on $Y=\{g_1=\cdots=g_N=0\}$ admits a Siu decomposition 
$$ 
[Y] =\sum\limits_{1\leq i_1<\cdots <i_{n-q}\leq N} 
\, T_{i_1,...,i_{n-q}} \wedge \bigwedge_{l=1}^{n-q} dg_{i_l}\,, 
$$
where the $T_{i_1,...,i_{n-q}}$ are $(0,n-q)$ currents which are regular holonomic because 
of Coleff-Herrera type 
(\cite{by3:gnus, meo:gnus, bjo:gnus}). One can then 
proceed in 
$$
{\bf U}:= \{(\zeta,\overline \zeta)\,:\, \zeta \in U\}\subset 
\C^{2n}
$$
with blocks of variables $(\zeta,\overline \zeta)$ and profit from the fact that 
formally $\partial_\zeta$ and $\partial_{\overline \zeta}$ can be considered as 
derivations respect to independent sets of variables. 

\noindent 
Consider then the function of one complex variable defined by
\begin{eqnarray}
\lambda \mapsto \Phi_{Y,f}(\lambda ):=\|f\|^{2\lambda}\,[Y]\, .
\end{eqnarray}
This function (which is a $'{\cal D}^{(n-q,n-q)} (U)$-current valued 
function) is well defined and holomorphic in 
$\{\lambda\in \C\, ;\,  {\rm Re}\, \lambda >0 \}$. Thanks to the Bernstein-Sato 
relations 
(2.1), it can be continued to the whole complex plane as a meromorphic function. 
 The poles of this meromorphic extension are among 
strictly negative rational numbers. Furthermore, there is a true pole at any 
point $\lambda =-k$, $k\in \N^*$. 

\vskip 2mm
\noindent
In fact, we will need a more precise result, where the construction 
of the meromorphic continuation of (2.2) is involved. What we need is 
formulated in the following proposition.   
\vskip 2mm

\begin{prop} 
Let $Y$ be an irreducible $q$-dimensional analytic subset of 
$U \subset \C^n$  
and $f_1,...,f_m$ $m$ functions holomorphic in $U$.  
For any $k\in \{1,...,m\}$ and for any ordered subset 
${\cal I}\subset \{1,...,m\}$ with cardinal $k\leq \min (q,m)$, the 
$'{\cal D}^{(n-q,n-q+k)}$-valued map   
$$
\lambda \mapsto \lambda c_k \|f\|^{2(\lambda-k-1)} \, [Y] \wedge 
\overline \partial \|f\|^2 \wedge \Big( \sum\limits_{l=1}^k (-1)^{l-1} \overline 
{f_{i_l}} \bigwedge_{{j=1}\atop {j\not=l}}^k \overline {df_{i_j}}\Big) 
$$
\(which is holomorphic in ${\rm Re}\, \lambda > k+1$\) can be continued as a meromorphic 
map to the whole complex plane, with no pole at $\lambda=0$. Its value at 
$\lambda=0$ defines a residual $(n-q,n-q+k)$-current which is supported by the 
analytic set 
$Y\cap \{f_1=\cdots=f_m=0\}=Y\cap V(f)$ and denoted as 
\begin{eqnarray}
\varphi \in {\cal D}^{(q,q-k)} 
\mapsto 
\RR \left[\matrix {[Y]\wedge (\cdot) \cr 
f_{i_1},...,f_{i_k} \cr 
f_1,....,f_m}\right] \, (\varphi)= 
\left[\matrix {[Y]\wedge \varphi \cr 
f_{i_1},...,f_{i_k} \cr 
f_1,....,f_m}\right]\,. 
\end{eqnarray}
\end{prop}
\noindent
{\bf Proof.} Assume that $Y$ is defined (in $U$) by the equations $g_1=\cdots=g_N=0$ 
and that $\nu$ is the multiplicity of the ideal ${\cal O}_{U,y}$ generated by 
$g_1,...,g_N$ at a generic point $y\in Y$. Let $d=n-q$. One can conclude from 
\cite{meo:gnus} that $[Y]$ coincides with the value at $\mu=0$ of the 
meromorphic $'{\cal D}^{(d,d)}(U)$-valued map $\Psi_g$  
$$
\mu 
\buildrel {\Psi_g}\over {\mapsto} { \mu  (d-1)! \over 
(2i\pi)^{d} \, \nu}\,  
\|g\|^{2\mu} \, 
\overline \partial \log \|g\|^2
\wedge \partial \log \|f\|^2 \wedge 
\sum\limits_{{j_1<\cdots <j_{d-1}}\atop 
{1\leq j_j\leq N}} \, 
\bigwedge_{l=1}^{d-1} 
\Big({\overline {\partial g_{j_l}} 
\wedge \partial g_{j_l} \over \|g\|^2} \Big)\,.    
$$
In fact, in the general situation where $(g_1,...,g_N)$ define 
a $q$-purely dimensional cycle ${\cal Z}$ 
(non necessarily irreducible) in $U$, the integration current 
(with multiplicities) on ${\cal Z}$ can be expressed as the value at $\lambda=0$ 
of some meromorphic $'{\cal D}^{(d,d)}(U)$-valued function which can be made explicit in terms 
of $g_1,...,g_N$ (see theorem 3.1 in \cite{by3:gnus} for a proof in the algebraic case).  
Let ${\cal I}\subset \{1,...,m\}$ with cardinal $k\leq \min (q,m)$ and, 
for ${\rm Re}\, \lambda > k+1$,   
$$
\Theta_{f,{\cal I}} (\lambda):= 
\lambda \|f\|^{2(\lambda-k-1)} \,  
\overline \partial \|f\|^2 \wedge \Big( \sum\limits_{l=1}^k (-1)^{l-1} \overline 
{f_{i_l}} \bigwedge_{{j=1}\atop {j\not=l}}^k \overline {df_{i_j}}\Big)\,.
$$ 

\noindent
In order to prove the proposition, we can localise the problem and assume that 
the origin belongs to $Y \cap V(f)$. 
As in our previous work (see for example \cite{by2:gnus}, pages 32-33, or 
\cite{by3:gnus}, page 208) we 
construct an analytic $n$-dimensional manifold
${\cal X}$, a neighborhood $V$ of $0$ in $U$, a
proper map $\pi: {\cal X}\rightarrow V$ which 
realizes a local isomorphism between
$V\setminus \{f_1\cdots f_m\, g_1 \cdots g_N=0\}$
and ${\cal X}\setminus \pi^{-1}(\{f_1\cdots f_m \, g_1\cdots g_N=0\})$, such
that in local coordinates on ${\cal X}$ (centered at a point $x$),
one has, in the corresponding
local chart ${\cal U}_x$ around $x$,
\begin{eqnarray*} 
f_j\circ \pi (t)&=& u_j(t)\, t_1^{\alpha_{j1}}\cdots t_n^{\alpha_{jn}}=
u_j(t)\, t^{\alpha_j},
\  j=1,\dots,m \\
g_k \circ \pi (t) &=& v_{k}(t)\,  t_j^{\beta_{k1}} \cdots t_n^{\beta_{kn}}=
v_k(t)\,  t^{\beta_k},\ k=1,...,N  
\end{eqnarray*}
where the $u_j$, $j=1,...,m$ and the $v_k$, $k=1,...,N$, 
are non vanishing holomorphic functions in ${\cal U}_x$, 
at least one of the monomials $t^{\alpha_j}$, $j=1,...,m$ divides all 
of them (we will denote this monomial as $t^\alpha$), and 
at least one of the monomials $t^{\beta_k}$, $k=1,...,N$ divides all 
of them (we will denote this monomial as $t^\beta$). 

\noindent
When $\varphi$ is a $(q,q-k)$-test form with support in 
$V$, one has, for ${\rm Re}\, \lambda >>0$, 
$$
\int_{V\,\cap \, Y} 
\Theta_{f,{\cal I}} (\lambda) \wedge \varphi = \Bigg[ \int_{V} 
\Psi_g(\mu) \wedge 
\Theta_{f,{\cal I}} (\lambda) 
\wedge \varphi \Bigg]_{\mu=0}
$$
(the right hand side being continued as a meromorphic function of $\mu$ 
which has no pole at $\mu=0$). 
For $\lambda$ fixed with ${\rm Re}\, \lambda >>0$, one can rewrite for 
${\rm Re}\, \mu >>0$ the integral 
$$
\int_{V} 
\Psi_g(\mu) \wedge 
\Theta_{f,{\cal I}} (\lambda) 
\wedge \varphi
$$
as a sum of integrals of the form 
\begin{eqnarray}
\int_{{\cal U}_x} \pi^*[\Psi_g] (\mu) \wedge 
\pi^* [\Theta_{f,{\cal I}} (\lambda)] \wedge \rho \pi^*(\varphi)\,,
\end{eqnarray} 
where $\rho$ is a test-function in ${\cal U}_x$ which corresponds to a partition of unity 
for $\pi^*({\rm Supp}\, \varphi)$. We know from lemma 1.1 and lemma 1.2 in 
\cite{by3:gnus} that 
$$
\Big[\pi^* [\Psi_g(\mu)]\Big]_{\mu=0}= \Big[ \Psi_{g\circ \pi} (\mu)\Big]_{\mu=0} 
$$
is a positive $\partial$ and $\overline \partial$-closed current $\theta_{{\cal U}_x}$ 
in ${\cal U}_x$, 
which implies that, as soon as ${\rm Re}\, \lambda >>0$,  
$$
\Bigg[\int_{{\cal U}_x} \pi^*[\Psi_g] (\mu) \wedge 
\pi^* [\Theta_{f,{\cal I}} (\lambda)] \wedge \rho \pi^*(\varphi)\Bigg]_{\mu=0}
= \int_{{\cal U}_x} 
\theta_{{\cal U}_x} \wedge \pi^* [\Theta_{f,{\cal I}} (\lambda)] \wedge \rho \pi^*(\varphi)\,.  
$$   
On the other hand, in ${\cal U}_x$ and for ${\rm Re}\, \lambda >>0$,  
a straightforward computation leads to   
$$
\pi^* \Big[\Theta_{f,{\cal I}}\Big] (\lambda)
=\lambda {a^{2\lambda} |t^{\alpha}|^{2\lambda} 
\over t^{k \alpha}} \Big( \vartheta + \varpi \wedge
{\overline{d t^{\alpha}} \over \overline {t^\alpha}}\Big)\,,
$$
where $\vartheta$ and $\varpi$ are smooth differential forms in ${\cal U}_x$ 
(with respective types $(0,k)$ and $(0,k-1)$) and 
$a$ is a strictly positive real analytic function in ${\cal U}_x$. 
It follows from Stokes's theorem that 
\begin{eqnarray}
\int_{{\cal U}_x} \pi^* [\Theta_{f,{\cal I}}] (\lambda) 
\wedge \theta_{{\cal U}_x} \wedge 
\rho \pi^*(\varphi) = 
\int_{{\cal U}_x} 
{|t^{\alpha}|^{2\lambda} \over 
t^ {k \alpha}} \, \theta_{{\cal U}_x} \wedge \xi_\varphi(\rho\,;\, t,\lambda)\,,
\end{eqnarray}
where $(t,\lambda) \mapsto \xi_{\varphi} (\rho\,;\, t,\lambda)$ 
is a $(n-q,n-q)$-differential form with smooth coefficients (in $t$) depending 
holomorphically in $\lambda$. 

\noindent
One can see also that, for ${\rm Re}\, \mu >> 0$,  
$$
\pi^* [\Psi_g] (\mu)= 
\mu \, b^{2\mu} \, |t^{\beta}|^{2\mu} 
\, \Big( {\overline {d t^\beta} \over \overline {t^\beta}} + \eta_{(0,1)} \Big) 
\wedge \Big({{d t^\beta} \over {t^\beta}} + \eta_{(1,0)} \Big) \wedge \upsilon \,,
$$
where $b$ is a strictly positive real analytic function in ${\cal U}_x$, 
$\eta_{(0,1)},\, \eta_{(1,0)},\, \upsilon$ are smooth differential forms in 
${\cal U}_x$ with respective 
types $(0,1)\,,\, (1,0)$ and $(d-1,d-1)$. This implies that, if 
$t_{i_1},...,t_{i_s}$ are the coordinates that appear in $t^\beta$\,,
$$
\theta_{{\cal U}_x} = \sum\limits_{l=1}^s \, [t_{i_l}=0]\,  \wedge \omega_{i_l} \,,
$$
where $\omega_{i_l}$ is a smooth $(d-1,d-1)$-form in ${\cal U}_x$ and 
$[t_i=0]$ denotes the integration current (without multiplicities) on $\{t_i=0\}$. 
Therefore, 
for ${\rm Re}\, \lambda >> 0$, 
$$
\int_{{\cal U}_x} \pi^* [\Theta_{f,{\cal I}}] (\lambda) 
\wedge \theta_{{\cal U}_x} \wedge 
\rho \pi^*(\varphi)=
\sum\limits_{{l=1}\atop 
{(t_{i_l}, t^\alpha)=1}}^s 
\int_{\{t_{i_l}=0\}\, \cap\,  {\cal U}_x} 
{|t^\alpha|^{2\lambda}
\over t^{k \alpha}} \omega_{i_l} (t) \wedge 
\xi_\varphi(\rho\,;\, t,\lambda).
$$
Such a function of $\lambda$ can be continued to a meromorphic function 
in the whole complex plane, with no pole at $\lambda=0$ 
(using Stokes's theorem). The assertion of the 
proposition follows, since for ${\rm Re}\, \lambda >>0$,  
$$
\int_{V} 
\Psi_g(\mu) \wedge 
\Theta_{f,{\cal I}} (\lambda) 
\wedge \varphi
$$
is a sum of integrals of the form (2.4). $\diamondsuit $
\vskip 2mm 
\noindent
Keeping the notation from above one has the following corollary.
\begin{corollary}
Under the conditions of proposition $2.1$, the residual current defined by $(2.3)$ has the 
following properties\\
$1)$ For any $h\in H(U)$ such that
\begin{eqnarray}
\forall K \subset \subset U\cap Y\,, \exists C_K>0\,,\ |h|\leq C_K \|f\|\ {on}\ K\,, 
\end{eqnarray}
one has 
$$
\RR \left[\matrix {h^k [Y]\wedge (\cdot ) \cr 
f_{i_1},...,f_{i_k} \cr
f_1,....,f_m}\right]\equiv 0 
$$
\\
$2)$ If $h\in H(U)$ and
\begin{eqnarray}
h(z)=0,\; \forall z\in Y\cap V(f)\,, 
\end{eqnarray}
then one has 
$$
\RR \left[\matrix {{\overline h }[Y]\wedge (\cdot ) \cr 
f_{i_1},...,f_{i_k} \cr
f_1,....,f_m}\right]\equiv 0 
$$ 
\end{corollary}    
\vskip 2mm
\noindent
{\bf Proof.} 
Let us now suppose that $h$ satisfies (2.6). If we do not 
perform integration by parts as in (2.5), we have, for ${\rm Re}\, \lambda >>0$,  
\begin{eqnarray*} 
&&\int_{{\cal U}_x} \pi^* [\Theta_{f,{\cal I}}] (\lambda) 
\wedge \theta_{{\cal U}_x} \wedge 
\rho \pi^*( h^k \varphi) \\
&&= 
\lambda 
\sum\limits_{{l=1}\atop 
{(t_{i_l}, t^\alpha)=1}}^s 
\int_{\{t_{i_l}=0\}\, \cap\, {\cal U}_x} 
{a^{2\lambda} |t^{\alpha}|^{2\lambda} 
\over t^{k \alpha}} \Big( \vartheta + \varpi \wedge
{\overline{d t^{\alpha}} \over \overline {t^\alpha}}\Big) 
\wedge \omega_{i_l} (t) \wedge \rho \pi^*(h^k \varphi)\,. 
\end{eqnarray*}     
Condition (2.6) implies that there exists some positive constant $\kappa$ 
such that, for any $l=1,...,s$ with $t_{i_l}$ coprime with $t^\alpha$,  
$$
|\pi^* h(t_1,...,\buildrel {i_l}\over{0},...,t_n)|\leq \kappa |t^\alpha|\,,\quad 
t\in {\rm Supp}\, \rho\,,
$$
which implies that $t^{k\alpha}$ divides $(\pi^* h^k)_{\{|t_{i_l}=0\}}$ on the 
support of $\rho$. This implies that for such $h$, 
$$
\Bigg[\int_{{\cal U}_x} \pi^* [\Theta_{f,{\cal I}}] (\lambda) 
\wedge \theta_{{\cal U}_x} \wedge 
\rho \pi^*( h^k \varphi)\Bigg]_{\lambda=0}=0\,,
$$
which gives the first assertion of the corollary 
since 
$$
\int_{V} 
\Psi_g(\mu) \wedge 
\Theta_{f,{\cal I}} (\lambda) 
\wedge \varphi
$$
as a sum of integrals of the form (2.4). \\
If $h$ vanishes 
on $Y\cap V(f)$, then, for any $l=1,...,s$ such that 
$t_{i_l}$ is coprime with $t^{\alpha}$, any coordinate which divides 
$t^\alpha$ also divides $(\pi^* h)_{|\{t_{i_l}=0\}}$ on the support of $\rho$.  
This implies that any expression of the form 
$$
\int_{\{t_{i_l}=0\} \cap {\cal U}_x} 
{a^{2\lambda} |t^{\alpha}|^{2\lambda} 
\over t^{k \alpha}} \Big( \vartheta + \varpi \wedge
{\overline{d t^{\alpha}} \over \overline {t^\alpha}}\Big) 
\wedge \omega_{i_l} (t) \wedge \rho \pi^*(\overline h \varphi)   
$$
has in fact no antiholomorphic singularity (therefore 
has a meromorphic extension which is polefree at the origin).
 It follows that for such $h$, one has again 
$$
\Bigg[\int_{{\cal U}_x} \pi^* [\Theta_{f,{\cal I}}] (\lambda) 
\wedge \theta_{{\cal U}_x} \wedge 
\rho \pi^*( \overline h \varphi)\Bigg]_{\lambda=0}=0\,,
$$
which proves the remaining assertion of the corollary since again 
$$
\int_{V} 
\Psi_g(\mu) \wedge 
\Theta_{f,{\cal I}} (\lambda) 
\wedge \varphi
$$
is a sum of integrals of the form (2.4). $\quad \diamondsuit$ 
\vskip 2mm
\noindent
When $k=m\leq q$, we will use the simplified notation 
$$
\RR \left[\matrix {[Y]\wedge (\cdot) \cr 
f_1,...,f_m }\right] \, (\varphi):= 
\RR \left[\matrix {[Y]\wedge (\cdot) \cr 
f_1,...,f_m\cr 
f_1,....,f_m}\right]\,.
$$
\vskip 2mm
\noindent
The transformation law for residual currents can be transposed to the 
case of restricted residual currents. Since we  deal in this paper with 
restricted residual currents supported by discrete sets, we state the 
transformation law in this particular setting. One has the following 
proposition~:
\vskip 2mm
\begin{prop}  
Let $Y$ be an irreducible $q$-dimensional analytic subset of 
$U \subset \C^n$ and $f_1,...,f_q,\, g_1,...,\,g_q$, $2q$ functions holomorphic in $U$ such that 
$Y\cap V(f)$ and $Y\cap V(g)$ are discrete analytic sets. Assume that there exist 
$q^2$ holomorphic functions in $U$, $a_{kl}$, $1\leq k,l\leq q$, such that 
$$
g_k(\zeta)=\sum\limits_{l=1}^q a_{kl} (\zeta) \, f_l(\zeta)\,
,\ k=1,...,q\,, \quad \zeta \in Y 
$$
Then, one has the following equality between restricted residual currents~: 
\begin{eqnarray} 
\RR \left[\matrix {[Y]\wedge (\cdot) \cr  
f_1,....,f_q}\right] = 
\RR \left[\matrix { \Delta \, [Y]\wedge (\cdot)\cr  
g_1,....,g_q}\right]\,,  
\end{eqnarray} 
where $\Delta := \det [a_{kl}]_{1\leq k,l\leq q}$. 
\end{prop} 

\noindent
{\bf Proof.} In order to prove this equality, we just need to prove it when $U$ is a neighborhood $V$  
of a point $\alpha\in Y\cap (V(f)\cup V(g))$ such that $\alpha$ is the only point of $Y\cap(V(f)\cup V(g))$ 
which lies in this neighborhhood. Thanks to the first assertion in Corollary 2.1,
 it is enough 
to test the two currents involved in (2.8) on test forms in 
${\cal D}^{(q,0)}(V)$ whose coefficients are holomorphic in a neighborhood of 
$\alpha$. Let $\varphi$ be such a test form. Since 
\begin{eqnarray*}
&&\overline \partial 
\Bigg[\|f\|^{2(\lambda-q)} \, [Y] \wedge  
\Big(\sum\limits_{j=1}^q  (-1)^{j-1} \overline f_j \bigwedge\limits_{{l=1}\atop {l\not=j}}^q 
\overline {df_j}\Big)\Bigg]\\
&&\qquad=\lambda \|f\|^{2(\lambda-q)} \, [Y]\, \wedge \bigwedge_{j=1}^q 
\overline {df_j} \\
&&\qquad= \lambda \|f\|^{2(\lambda-q-1)}\, [Y] 
\wedge \overline\partial \|f\|^2 
\wedge  
\Big(\sum\limits_{j=1}^q  (-1)^{j-1} \overline f_j \bigwedge\limits_{{l=1}\atop {l\not=j}}^q 
\overline {df_l}\Big)
\end{eqnarray*}
for ${\rm Re}\, \lambda >>0$ and
$$
\sum\limits_{j=1}^q s_j (\zeta) f_j(\zeta)=1\,,\quad \forall\, 
\zeta \in (V\cap Y_{\rm reg})\setminus \{\alpha\}, 
$$ 
where 
$$
s_j:= {\overline f_j \over \|f\|^2}\,,\quad j=1,...,q, 
$$
one has, by Stokes's theorem, that 
\begin{eqnarray} 
\RR \left[\matrix {[Y]\wedge \varphi \cr  
f_1,....,f_q}\right]&=& 
(-1)^q \omega_q 
\int_{Y_{\rm reg}} {\sum\limits_{j=1}^q  (-1)^{j-1} \overline f_j \bigwedge\limits_{{l=1}\atop {l\not=j}}^q 
\overline {df_j}\over \|f\|^{2q}} \wedge \overline\partial \varphi \nonumber \\
&=& (-1)^q \omega_q 
\int_{Y_{\rm reg}} \Big(\sum\limits_{j=1}^q  (-1)^{j-1} s_j\bigwedge\limits_{{l=1}\atop {l\not=j}}^q 
{ds_l}\Big)\wedge \overline\partial \varphi \, .
\end{eqnarray}
Similarly, if we introduce 
$$
t_j:={\overline g_j \over \|g\|^2}\,,\quad j=1,...,q\,,
$$
and 
$$
\widetilde {s_j}:=\sum\limits_{l=1}^q a_{lj} t_l\,,\quad j=1,...,q\,,
$$
one has also 
$$
\sum\limits_{j=1}^q \widetilde {s_j} (\zeta) f_j(\zeta)=1\,,\quad \forall\, 
\zeta\in (V \cap Y_{\rm reg})\setminus \{\alpha\}\, .  
$$ 
Let, for $\xi\in [0,1]$ and $j=1,...,q$, 
$$
s_j^{(\xi)} = (1-\xi)\, s_j + \xi\,  \widetilde  {s_j}\, . 
$$
Note that we have 
$$
\sum\limits_{j=1}^q s_j^{(\xi)} (\zeta) f_j(\zeta) =1\,,\quad 
\forall \xi \in [0,1]\,,\quad 
\forall \zeta 
\in (V \cap Y_{\rm reg})\setminus \{\alpha\}\, . 
$$ 
Therefore, one has, since 
$$
\bigwedge\limits_{j=1}^q \overline\partial_\zeta  s_j^{(\xi)}\equiv 0 
$$
on $(V \cap Y_{\rm reg})\setminus \{\alpha\}$\,,    
$$
{d\over d\xi} \Bigg[ 
\int_{W_{\rm reg}} \Big(\sum\limits_{j=1}^q  (-1)^{j-1} s_j^{(\xi)}\bigwedge\limits_{{l=1}\atop {l\not=j}}^q 
{ds_l^{(\xi)}}\Big)\wedge \overline\partial \varphi \Bigg] \equiv 0 
$$
on $[0,1]$. It follows from (2.9) that 
\begin{eqnarray*}
\RR \left[\matrix {[Y]\wedge \varphi \cr  
f_1,....,f_q}\right] &=& 
(-1)^q \omega_q \int_{Y_{\rm reg}} \Big(\sum\limits_{j=1}^q  (-1)^{j-1} \widetilde {s_j}\bigwedge\limits_{{l=1}\atop {l\not=j}}^q 
{\overline \partial \widetilde {s_l}}\Big)\wedge \overline\partial \varphi \\
&=& (-1)^q \omega_q 
\int_{Y_{\rm reg}} \Delta\, {\sum\limits_{j=1}^q  (-1)^{j-1} \overline g_j \bigwedge\limits_{{l=1}\atop {l\not=j}}^q 
\overline {dg_j}\over \|g\|^{2q}} \wedge \overline\partial \varphi \\
&=& \RR \left[\matrix {\Delta\, [Y]\wedge \varphi \cr  
g_1,....,g_q}\right] \, . 
\end{eqnarray*}
this concludes the proof of the 
proposition. $\quad\diamondsuit$ 
\vskip 2mm
\noindent
As a consequence of this result, we will state in the algebraic context the following 
analog of the global transformation law. We need first some piece of notation.
 Assume that 
$W$ is a $q$-dimensional irreducible algebraic subvariety in the affine space 
${\bf A}^n_\C$ 
(the integration current on $|W|$ without multiplicities taken into account being
 denoted 
as $[W]$) and that $P_1,...,P_q$ are $q$ elements in $\C[X_1,...,X_n]$ such that 
$|W|\cap V(P_1,...,P_q)$ is a discrete (hence finite) algebraic subset of 
$\C^n$. For 
any $Q\in \C[X_1,...,X_n]$, any ordered subset $\{i_1,...,i_q\}$ of 
$\{1,...,n\}$, we will denote as 
$$
\RR \left[\matrix{ [W]\wedge Q \bigwedge\limits_{l=1}^q dX_{i_l} \cr \cr 
P_1,...,P_q}\right] 
$$
the result of the action of the $W$-restricted current 
$$
\varphi \mapsto  \RR \left[\matrix{ [W]\wedge \varphi \cr 
P_1,...,P_q}\right]
$$
on the $(q,0)$-test form $Q(\zeta) \psi(\zeta) \bigwedge_{l=1}^q d\zeta_{i_l}$,
 where $\psi$ is any 
test-function in ${\cal D} (\C^n)$ which equals $1$ in a neighborhood of $|W|\cap V(P)$. If 
$$
\Gamma = \sum\limits_{j=1}^M \nu_j W_j 
$$
(where $W_1,...,W_M$ are 
$M$ irreducible algebraic subsets in $\C^n$ and $\nu_j\in \N^*$, $j=1,...,M$)  
is an effective $q$-dimensional algebraic cycle in the affine space $\C^n$ and
 $P_1,...,P_q$ are 
$q$ polynomials such that $W_j\cap V(P)$ is discrete for any $j=1,...,M$, we will
 also denote as  
$$
\RR \left[\matrix{ [\Gamma]\wedge Q \bigwedge\limits_{l=1}^q dX_{i_l}  \cr\cr 
P_1,...,P_q}\right]
$$
the weighted sum 
$$
\sum\limits_{j=1}^M \nu_j 
\RR \left[\matrix{ [W_j]\wedge Q \bigwedge\limits_{l=1}^q dX_{i_l} \cr \cr
P_1,...,P_q}\right]
$$
\vskip 2mm
\begin{corollary}
Let $\Gamma$ be an effective $q$-dimensional algebraic cycle in the affine 
space $\C^n$ and 
$P_1,...,P_q, R_1,...,R_q$ be $2q$ polynomials such that 
${\rm Supp}\, \Gamma \cap 
V(P_1,...,P_q)$ and ${\rm Supp}\, \Gamma \cap V(R_1,...,R_q)$ 
are discrete (hence finite) 
algebraic subsets of $\C^n$. Assume that there is a $(q,q)$-matrix of polynomials 
$[A_{k,l}]_{1\leq k,l\leq q}$ such that 
$$
R_k =\sum\limits_{l=1}^q A_{kl} P_l \quad on\quad {\rm Supp}\,
 \Gamma\,\quad k=1,..,q\,.
$$
Then, for any $Q\in \C[X_1,...,X_n]$, any ordered subset $\{i_1,...,i_q\}$ of 
$\{1,...,n\}$, one has 
\begin{eqnarray}
\RR \left[\matrix{ [\Gamma]\wedge Q \bigwedge\limits_{l=1}^q dX_{i_l} \cr \cr 
P_1,...,P_q}\right] =
\RR \left[\matrix{ [\Gamma]\wedge \Delta \, Q \bigwedge\limits_{l=1}^q dX_{i_l}  \cr\cr  
R_1,...,R_q}\right]\,,
\end{eqnarray}
where $\Delta$ denotes the determinant of the matrix $[A_{k,l}]_{1\leq k,l\leq q}$. 
\end{corollary}  
\vskip 2mm
\noindent
Another key point about the restricted residual current in the discrete context 
is the following annihilating property~: 
\vskip 2mm
\begin{prop}  
Let $Y$ be an irreducible $q$-dimensional analytic subset of 
$U \subset \C^n$ and $f_1,...,f_q$ be $q$ functions holomorphic in $U$ such that 
$Y\cap V(f)$ is a discrete analytic set. Then 
 one has, for $k=1,...,q$, 
\begin{eqnarray} 
{\rm Res}\, \left[\matrix {f_k[Y]\wedge (\cdot) \cr  
f_1,....,f_q }\right] = 0
\end{eqnarray} 
\end{prop}

\noindent
{\bf Proof.} 
We give here a self-contained proof of the above proposition.
 Actually, because 
of the properties quoted in Corollary 2.1, it is 
enough to show that if $\alpha\in V(P)\cap Y$ and $\varphi$ is a test-function    
with support arbitrarly small about $\alpha$ with $\varphi=1$ in some neighborhhood 
$v_\alpha$ of $\alpha$, then, for any function $h \in C^\infty (U)$ which is 
holomorphic on $v_\alpha$, for any ordered subset ${\cal I}=\{i_1,...,i_q\}\subset \{1,...,n\}$, 
one has, for $j=1,...,q$, 
$$
{\rm Res}\, \left[\matrix {f_j[Y]\wedge h \varphi \, d\zeta_{\cal I} \cr  
f_1,....,f_q}\right]=0\,. 
$$ 
One can use Stokes's formula (as in the proof of proposition 2.2) and write 
$$
{\rm Res}\, \left[\matrix {f_k[Y]\wedge h \varphi \, d\zeta_{\cal I} \cr  
f_1,....,f_q}\right]=
(-1)^q \omega_q \int_{Y_{\rm reg}} h f_k  
\Big( \sum\limits_{j=1}^q (-1)^{j-1} s_j \bigwedge\limits_{{l=1}\atop {l\not=j}}^q ds_l \Big)
\wedge \overline \partial \varphi \wedge d\zeta_{\cal I}
\,,
$$ 
where $s_j:=\overline {f_j}/\|f\|^2$, $j=1,...,q$. 
One can see at once that 
\begin{eqnarray*}
&&f_k 
\Big( \sum\limits_{j=1}^q (-1)^{j-1} s_j \bigwedge\limits_{{l=1}\atop {l\not=j}}^q ds_l \Big)
\wedge d\zeta_{\cal I} \wedge [Y]= 
\Big(\bigwedge_{{l=1}\atop {l\not=k}}^q ds_l \Big) \wedge 
\overline \partial \varphi \wedge d\zeta_{{\cal I}} \wedge [Y]\\
&&\qquad\qquad= \pm d\Bigg[ s_{k'} 
\Big(\bigwedge_{{l=1}\atop {l\not=k,k'}}^q ds_l \Big)
\wedge 
\overline \partial \varphi \wedge d\zeta_{{\cal I}} \wedge [Y]\Bigg] 
\end{eqnarray*}
for $k'\not=k$, 
since $s_1 f_1+\cdots +s_q f_q\equiv 1$ on 
$Y\cap {\rm Supp}\ \overline\partial\varphi$, 
which shows that 
$$
\int_{Y_{\rm reg}} h f_k  
\Big( \sum\limits_{j=1}^q (-1)^{j-1} s_j \bigwedge\limits_{{l=1}\atop {l\not=j}}^q ds_l \Big)
\wedge \overline \partial \varphi \wedge d\zeta_{\cal I} =0
$$
as a consequence of Stokes's formula on $Y$. $\quad \diamondsuit$
\vskip 2mm
\noindent
We remark here that there is an alternative proof of the last 
proposition. In fact,  when $m\leq q$ and $f_1,...,f_m$ define a complete intersection 
on $Y$, one can show that the restricted residual current 
$$
{\rm Res}\, \left[\matrix{ [Y] \wedge (\cdot) \cr f_1,...,f_m}\right] 
$$
coincides with the Coleff-Herrera current 
$\Big(\bigwedge_{j=1}^m \overline\partial (1/f_j)\Big) 
\wedge [Y]$ 
as it is defined in \cite{coh:gnus}. The proof of this claim  can be carried out
 as it is done 
in the non restricted case in \cite{pty:gnus}, section 4. Since the proof of this 
fact is rather tedious, we will not give it here. 
A consequence of this result is that, when 
$f_1,...,f_m$ ($m\leq q$) define a complete intersection on $Y$, one has, 
for $k=1,...,m$,  
$$
{\rm Res}\, \left[\matrix {f_k[Y]\wedge (\cdot) \cr  
f_1,....,f_m}\right] = f_k \Big(\bigwedge_{j=1}^m \overline\partial {1\over f_j} \Big)\wedge [Y]=0
$$
(see \cite{coh:gnus}). This implies the proposition when $m=q$. 

\noindent

\section{An Abel-Jacobi formula in the restricted case (analytic approach)} 
\setcounter{equation}{0}

One of the key facts about restricted residual currents (as defined through the 
analytic approach described in section $2$) is that they satisfy (in the $0$-dimensional complete 
intersection setting) Abel-Jacobi's formula, exactly as in the non-restricted case
(see \cite{vy:gnus}). Such a result 
will be, together with the validity of the transformation law in the restricted context) 
a crucial fact in order to compare our analytic approach and the algebraic one.  

\begin{prop} 
Let $W$ be a $q$-dimensional irreducible affine algebraic subvariety of 
the affine scheme ${\bf A}^n_\C$ 
\($0<q<n$\) 
and $P_1,...,P_q$ be $q$ polynomials in $\C[X_1,...,X_n]$ such that there exists 
strictly positive rational numbers $\delta_1,...,\delta_q$ and two constants $K>0$,
 $\kappa>0$ 
with~: 
\begin{eqnarray} 
\zeta \in |W|\,,\ \|\zeta\|\geq K \ \Longrightarrow \ \sum\limits_{j=1}^q {|P_j(\zeta)|\over 
\|\zeta\|^{\delta_j}} \geq \kappa\,.
\end{eqnarray} 
Then, for any $Q\in \C[X_1,...,X_n]$ such that ${\rm deg}\, Q <\delta_1+\cdots +\delta_q -q$,
for any multi-index $(i_1,...,i_q)$ in $\{1,...,n\}^q$,  
\begin{eqnarray}
\RR \left[\matrix{ [W]\wedge Q \bigwedge\limits_{l=1}^q dX_{i_l}\cr  \cr 
P_1,...,P_q}\right] =0\, .
\end{eqnarray}
\end{prop} 

\noindent
Before we give the proof of this result, let us state an important corollary~:  
\begin{corollary} 
Let $W$ be a $q$-dimensional irreducible algebraic subvariety in the affine scheme 
${\bf A}^n_\C$ 
and ${\cal W}$ be its completion in ${\rm Proj}\, \C[X_0,...,X_n]$. 
Let $P_1,...,P_q$ be $q$ 
elements in $\C[X_1,...,X_n]$, with respective degrees $D_1,...,D_q$, such that 
\begin{eqnarray}
&&|{\cal W}| \cap \Big\{[\zeta_0\,:\, ...\,:\, \zeta_n]\in \P^n(\C)\,;\, 
{}^h P_j (\zeta_0,...,\zeta_n)=0,\; j=1,...,q \Big\} \subset \C^n \,,\nonumber \\
&& 
\end{eqnarray}
where ${}^h P_j$, $j=1,...,q$, denotes the homogeneization of the 
polynomial $P_j$. Then, for any $Q\in \C[X_1,...,X_n]$ such that ${\rm deg}\, Q <D_1+\cdots +D_q -q$,
for any multi-index $(i_1,...,i_q)$ in $\{1,...,n\}^q$,  
\begin{eqnarray}
\RR \left[\matrix{ [W]\wedge Q \bigwedge\limits_{l=1}^q dX_{i_l}\cr  \cr 
P_1,...,P_q}\right] =0
\end{eqnarray} 
\end{corollary} 

\vskip 2mm
\noindent
{\bf Proof of corollary 3.1.} Assume that 
$$
|{\cal W}|=\Big\{[\zeta_0\,:\, ...\,:\, \zeta_n]\in \P^n(\C)\,;\, 
{}^h {\cal G}_j (\zeta_0,...,\zeta_n)=0,\; j=1,...,N \Big\}\,,
$$
where ${\cal G}_1,...,{\cal G}_N$ are homogeneous polynomials in $\widetilde \zeta= (\zeta_0,...,\zeta_n)$. 
Condition (3.3) implies that 
$$
|{\cal W}| \cap \Big\{[\zeta_0\,:\, ...\,:\, \zeta_n]\in \P^n(\C)\,;\, 
{}^h P_j (\zeta_0,...,\zeta_n)=0,\; j=1,...,q \Big\} 
$$
is a finite set in $\C^n$~; this implies (through a compacity argument) 
that there exists $K, \kappa>0$ such that, for any 
$(\zeta_0,...,\zeta_n) \in \C^{n+1}\setminus \{(0,...,0)\}$ such that 
$$
\Big(|\zeta_1|^2 +\cdots + |\zeta_n|^2\Big)^{1/2} \geq K |\zeta_0|\,,
$$
one has 
$$
\sum\limits_{j=1}^q {|{}^h P_j(\widetilde \zeta)|\over 
\|\widetilde \zeta \|^{D_j}} 
+ \sum\limits_{l=1}^M {|{\cal G}_l(\widetilde \zeta)|\over 
\|\widetilde \zeta\|^{{\rm deg}\, {\cal G}_l} } \geq \kappa \,.
$$
Condition (3.1) with 
$\delta_j=D_j$, $j=1,...,q$, holds if we restrict to the affine space $\C^n$.   
The statement (3.4) follows then from (3.2). $\quad \diamondsuit$ 

\vskip 2mm
\noindent
We remark that a proposition similar to proposition 3.1 was proved in the 
non restricted case ($W={\bf A}^n_\C$) in \cite{vy:gnus}. Unfortunately, the proof 
which is given there 
(and depends heavily on resolution of singularities on the analytic manifold 
$\P^n(\C)$) 
cannot immediately be transposed to the restricted case 
(since the Zariski closure $|{\cal W}|$ of $|W|$ in $\P^n(\C)$ is not a smooth 
manifold anymore). Instead, we will 
follow an alternative approach (applicable also for the case $q=n$), 
based on an argument in the affine space 
(and not in its compactification $\P^n(\C)$), which was proposed by 
Ha{\"\i} Zhang in \cite{zhang:gnus}. Our task has been to adapt this argument to 
the restricted case. 
\vskip 2mm
\noindent
Note that, if $z=Aw$ is a linear change of variables in $\C^n$, one has, for any element  
in ${\cal D}^{(q,0)}(\C^n)$
$$
\RR
\left[\matrix { [W]\wedge \varphi \cr f_1,...,f_q}\right]= 
\RR 
\left[ \matrix{ 
[A^{-1} (W)] \wedge A^* \varphi \cr f_1\circ A,...,f_q\circ A}\right]\, .
$$  
Therefore, we do not loose generality is we assume that ${\cal I}=\{1,...,q\}$ and that the projection 
$$
\Pi~: (\zeta_1,...,\zeta_n) \mapsto (\zeta_1,...,\zeta_q) 
$$
is a proper map from $|W|$ to $\C^q$ (coordinates can be choosen in such a way that 
Noether normalization theorem applies respect to any $(q,n-q)$ splitting 
$\zeta=(\zeta',\zeta'')$ of the set of variables $(\zeta_1,...,\zeta_n)$, see 
for example \cite{fl:gnus,rud:gnus}).       
\vskip 2mm
\noindent    
For $\delta _i$, $i=1,\dots , q $ which appear in the 
statement of Proposition 3.1 we choose a positive integer $N$ large enough so that 
\begin{eqnarray}
N\prod\limits _{{l=1}\atop {l\not=j}}^q\delta_l >2, \quad 
j=1, \dots , q \, . 
\end{eqnarray} 
Then, let 
$$
\delta^{[j]}:=
N\prod\limits _{{l=1}\atop {l\not=j}}^q \delta_l\,,\quad j=1,...,q\,,
$$
and 
$$
\delta:= N \delta_1 \cdots \delta_q=\delta_j \delta^{[j]}\,,\quad j=1,...,q\, .
$$
Similarly, for the polynomials 
$P_1, \dots, P_q $, one can define, in the affine open set 
$\C^n \setminus \{P_1\cdots P_q=0\}$, the 
$C^\infty$ functions 
$$
\widetilde s_j:= {|P_j|^{\delta^{[j]}} 
\over P_j\, \sum\limits_{l=1}^q |P_l|^{\delta^{[l]}}}\,,\quad j=1,...,q\, .
$$
These functions $\widetilde s_j$, $j=1,...,q$, 
extend (provided $N >> 1$) to $C^1$ functions in 
$\C^n \setminus V(P)$, satisfying 
$$
\sum\limits_{j=1}^q \widetilde s_j (\zeta) P_j(\zeta)=1\,,\quad 
\zeta\in \C^n \setminus V(P)\, . 
$$
Let finally 
$$
u_j:=|P_j|^{\delta^{[j]}/2}\,,\quad j=1,...,q
$$
and 
$$
S:=\sum\limits_{j=1}^q u_j^2 =\| u\|^2\,.
$$
\vskip 2mm
\noindent
At this point we return to the \\
\\
\noindent
{\bf Proof of proposition 3.1.} 
One can suppose without any loss of generality that $\{i_1,...,i_q\}=
\{1,...,q\}$ and that the projection $\Pi$ is a proper map from 
$|W|$ to $\C^q$. 
Condition (3.1) implies the existence of a strictly positive constant 
$\kappa_N$ such that
\begin{eqnarray}
S(\zeta) \geq \kappa_N \|\zeta\|^{\delta}\,,\quad 
\zeta \in |W|\,,\ \|\zeta\|\geq K\, .
\end{eqnarray}
\vskip 2mm
\noindent
Let 
$$
\theta \in {\cal D}(]-3\kappa_N/4,3\kappa_N/4[
$$ 
such that $\theta\equiv 1$ on $[-\kappa_N/4,\kappa_N/4]$~; for any $R>0$, 
let the element $\varphi_R$ in 
$C^1(\C^n)$ defined as  
$$
\varphi_R~: \zeta \mapsto \theta (S(\zeta) /R^\delta)\,.
$$
Since the restriction $S_{|_{|W|}}$ is a proper map 
(all $\delta_j$'s, $j=1,...,q$, being strictly positive) and 
$V(P)\cap |W|$ is a discrete (hence finite) 
algebraic subset of $\C^n$ (this follows also from (3.1)), there exists $R_0$ such that 
for $R>R_0$, $\varphi_R \equiv 1$ in a neighborhood of $|W|\cap V(P)$. Therefore, if  
$$
s_j:={\overline {P_j} \over \|P\|^2} 
$$
one has (see for example formula (2.9))
\begin{eqnarray}
&&\RR \left[\matrix{ [W]\wedge Q \bigwedge\limits_{l=1}^q dX_{l}\cr  \cr 
P_1,...,P_q}\right]\nonumber \\
&&\qquad\qquad\qquad = c_q    
\int_{|W|_{\rm reg}} \Big(\sum\limits_{j=1}^q  (-1)^{j-1} s_j\bigwedge\limits_{{l=1}\atop {l\not=j}}^q 
{ds_l}\Big)\wedge Q\, d\zeta'\, \wedge \overline\partial\varphi_R
\nonumber \\  
\end{eqnarray} 
for any $R>R_0$, where $d\zeta'=\bigwedge_{l=1}^q d\zeta_l$. 
It follows from an homotopy argument similar to the one  which is developped in the 
proof of proposition 2.2 that 
\begin{eqnarray}
&&\RR \left[\matrix{ [W]\wedge Q \bigwedge\limits_{l=1}^q dX_{l}\cr  \cr 
P_1,...,P_q}\right]\nonumber \\
&&\qquad\qquad\qquad =c_q   
\int_{|W|_{\rm reg}} \Big(\sum\limits_{j=1}^q  (-1)^{j-1} \widetilde 
s_j\bigwedge\limits_{{l=1}\atop {l\not=j}}^q 
{d\widetilde 
s_l}\Big)\wedge Q\, d\zeta'\, \wedge 
\overline\partial\varphi_R
\nonumber \\
&&\qquad\qquad\qquad =c_q    
\int_{|W|_{\rm reg}} \Big(\sum\limits_{j=1}^q  (-1)^{j-1} \widetilde 
s_j\bigwedge\limits_{{l=1}\atop {l\not=j}}^q 
{\overline\partial \widetilde 
s_l}\Big)\wedge Q \, d\zeta'\, \wedge 
\overline\partial\varphi_R\nonumber \\
\end{eqnarray} 
for any $R>R_0$. 
Since $ P_j \widetilde s_j= u_j^2/S$, $j=1,...,q$, one can rewrite (3.8) as
\begin{eqnarray}
&&\RR \left[\matrix{ [W]\wedge Q \bigwedge\limits_{l=1}^q dX_{l}\cr  \cr 
P_1,...,P_q}\right]\nonumber \\
&&=c_q\,  2^{q-1}    
\int_{|W|_{\rm reg}}
\Big(\prod\limits_{j=1}^q {|P_j|\over P_j}\, 
u_j^{1-{2\over \delta^{[j]}}}\Big) 
{\sum\limits_{j=1}^q (-1)^{j-1} 
u_j\bigwedge\limits_{{l=1}\atop {l\not=j}}^q 
{d 
u_l} \over \|u\|^{2q}} \wedge Q\, d\zeta'\, \wedge 
\overline\partial\varphi_R 
\nonumber \\
&&={ (-1)^q \, c_q\, 2^{q-1}\over R^\delta} 
\int_{|W|_{\rm reg}}
\Big(\prod\limits_{j=1}^q {|P_j|\over P_j}\, 
u_j^{1-{2\over \delta^{[j]}}}\Big) 
\, {\bigwedge\limits_{l=1}^q du_l \over \|u\|
^{2(q-1)}} 
\wedge \theta'\Bigg({\|u\|^2\over R^\delta}\Bigg) \, 
Q\, d\zeta'. \nonumber \\
\end{eqnarray}
\vskip 2mm
\noindent
For any order ${\cal J} \subset \{1,...,q\}$, let 
$$
\omega_{\cal J}=\bigwedge\limits_{j=1}^q d\alpha_{{\cal J},l}
$$
where 
$$
\alpha_{{\cal J},l} (\zeta_1,...,\zeta_q):= 
\cases {{\rm Re}\, \zeta_j\ {\rm if}\ j\in {\cal J} \cr
{\rm Im}\, \zeta_j\ {\rm if}\, j\not\in {\cal J}}\,;
$$
then one can write 
$$
d\zeta'=\bigwedge\limits_{l=1}^q d\zeta_l =\sum\limits_{{\cal J}\subset \{1,...,q\}}
i^{q-\#{\cal J}}\, d\omega_{{\cal J}}\, . 
$$   
In order to prove formula (3.4), it is enough to prove that for any ${\cal J}
\subset \{1,...,q\}$, one has, as soon as ${\rm deg}\, Q <\delta_1+\cdots+\delta_q -q$,  
\begin{eqnarray} 
&&\lim\limits_{R\rightarrow +\infty} 
\Bigg[{1\over R^\delta}   
\int_{|W|_{\rm reg}}
\Big(\prod\limits_{j=1}^q {|P_j|\over P_j}\, 
u_j^{1-{2\over \delta^{[j]}}}\Big) 
\, {\bigwedge\limits_{l=1}^q du_l \over \|u\|
^{2(q-1)}} 
\wedge \theta'\Bigg({\|u\|^2 \over R^\delta}\Bigg) \, 
Q\, \omega_{{\cal J}} (\zeta') \Bigg] =0\, .\nonumber\\
&&
\end{eqnarray}
Since the restriction of $P=(P_1,...,P_q)$ to each connected sheet 
${\cal F}$ (above the $\zeta'$-space) of the $2q$-dimensional  
real manifold 
$|W|_{\rm reg}$ is proper, the map 
$$
F_{\cal J,{\cal F}}~:\ 
\zeta \in {\cal F}\mapsto (u_1,...,u_q, \alpha_{{\cal J},1},..., 
\alpha_{{\cal J},q})  
$$
is a $\R^{2q}$-valued proper map, with topological degree $d_{{\cal J},{\cal F}}$. 
Moreover, condition (3.6) implies that, for $R>K$, 
$$
{\rm Supp} \Big(\theta(S/R^{\delta})\Big)
\subset \{ \zeta \in \C^n\,:\, \|\zeta\| <R\}\, . 
$$
 Actually, for $\|\zeta\|\geq R>K$, one has 
 $$
 S(\zeta) \geq \kappa_N \|\zeta\|^{\delta} 
\geq \kappa_N R^{\delta} > (3 \kappa_N/4) R^{\delta}\, .
$$
 For such $R$, one has 
$$
\Bigg\|\prod\limits_{j=1}^q {|P_j|\over P_j} \, Q\, \theta' (S/R^{\delta})\Bigg\|_\infty 
\leq C\, R^{\,{\rm deg}\, Q}\,,
$$
where $C=C(\theta, Q)$ is a positive constant. It follows then from the properness of 
all maps   
$F_{{\cal J},{\cal F}}$ and from the positivity 
of the differential form 
$$
\Big(\prod\limits_{j=1}^q \, 
u_j^{1-{2\over \delta^{[j]}}}\Big) 
\bigwedge\limits_{l=1}^q du_l
$$
in $]0,\infty[^q$ that 
\begin{eqnarray*}
&&{1\over R^\delta} \Bigg| 
\int_{|W|_{\rm reg}}
\Big(\prod\limits_{j=1}^q {|P_j|\over P_j}\, 
u_j^{1-{2\over \delta^{[j]}}}\Big) 
\, {\bigwedge\limits_{l=1}^q du_l \over \|u\|
^{2(q-1)}} 
\wedge \theta'\Bigg({\|u\|^2 \over R^\delta}\Bigg) \, 
Q\, \omega_{{\cal J}} (\zeta') \Bigg|\\  
&& \quad 
\leq 
{(\sum\limits_{\cal F} d_{{\cal J},{\cal F}})\, C \, R^{\, {\rm deg}\, Q} \over R^{\delta}} 
\Bigg( \int_{{\kappa_N R^{\delta}\over 4} 
\leq \|u\|^2 \leq {3 \kappa_N R^{\delta}\over 4}} 
\Big( \prod\limits_{j=1}^q u_j^{1-{2\over \delta^{[j]}}}\Big) 
{\bigwedge\limits_{l=1}^q du_l 
\over \|u\|^{2(q-1)}}
 \Bigg) \\
&&\quad\quad \times \Bigg(
\int_{\|t\| < R} dt_1\wedge \cdots \wedge dt_q \Bigg) \\
&&\leq 
{(\sum\limits_{\cal F} d_{{\cal J},{\cal F}})\,  \, C_N \, R^{\, {\rm deg}\, Q + q} \over R^{q \delta}}   
\Bigg( \int_{{\kappa_N R^\delta \over 4} 
\leq \|u\|^2 \leq {3 \kappa_N R^\delta \over 4}} 
\Big( \prod\limits_{j=1}^q u_j^{1-{2\over \delta^{[j]}}}\Big) 
\bigwedge\limits_{l=1}^q du_l \Bigg)\\
&&\leq 
{(\sum\limits_{\cal F} d_{{\cal J},{\cal F}})\,  \, \widetilde C_{N,\vec \delta} 
\, R^{\, {\rm deg}\, Q + q} \over R^{q \delta}} \, R^{
{\delta\over 2} \sum\limits_{j=1}^q \Big(1-{1\over \delta^{[j]}}\Big)+ 
q{\delta \over 2}}\\
&&\leq 
(\sum\limits_{\cal F} d_{{\cal J},{\cal F}})\,  \, \widetilde C_{N,\vec \delta}
\, R^{{\rm deg}\, Q +q -\delta_1-\,\cdots\, -\delta_q}= {\bf o} (1)\,,
\end{eqnarray*}   
which proves the conclusion (3.10) we need. The proof of proposition 3.1 
is therefore completed. $\quad\diamondsuit$ 

\section{Analytic versus algebraic approach} 
\setcounter{equation}{0}
Let ${\cal X}$ be an integral $\C$-variety of dimension $q$ and ${\cal D}_1$,...,${\cal D}_q$ be  
$q$ Cartier divisors on ${\cal X}$ such that $|{\cal D}_1|\cap \cdots \cap 
|{\cal D}_q|$ is finite. If $\omega$ is a meromorphic form 
in $\Omega^q_{\C({\cal X})/\C}$ which has a simple pole along  
${\cal D}_1+\cdots+ {\cal D}_q$,
 one may define (see \cite{hu:gnus}, page 621) the local    
residue of $\omega$ at any closed point 
$\alpha$ in $|{\cal D}_1|\cap \cdots \cap |{\cal D}_q|$. That is, 
if 
$$
\omega = {\eta \over f_1 \cdots f_q}\,,
$$
where $\eta \in \omega^q_{\C({\cal X})/\C,\alpha}$ and $f_j=0$, $j=1,...,q$, 
is a local equation for ${\cal D}_j$ at $\alpha$ then   
$$
{\rm Res}_{{\cal X};{\cal D}_1,...,{\cal D}_q,\alpha} \, (\omega) 
={\rm Res}_{\C({\cal X})/\C,\alpha} \, \Bigg( \Bigg[ \matrix{\eta \cr 
f_1,...,f_q }\Bigg]\Bigg)\,. 
$$
 When ${\cal X}$ is smooth, this definition agrees 
with the definition in \cite{gh:gnus}, chapter 5, section 1 (see \cite{lip1:gnus}, 
Appendix A). Adding the hypothesis that ${\cal X}$ is $\C$-complete, one has  
(see proposition 12.2, page 108, in \cite{lip2:gnus}) 
$$
\sum\limits_{\alpha \in |{\cal D}_1|\cap \cdots \cap |{\cal D}_q|}
\, {\rm Res}_{{\cal X};{\cal D}_1,...,{\cal D}_q,\alpha} \, (\omega)=0\,,
$$
which is known as residue theorem on ${\cal X}$ (it extends the classical 
residue theorem on a complete integral curve in its algebraic 
formulation, see \cite{se:gnus}). 
\vskip 2mm
\noindent
Such a residue theorem holds in our analytic setting (and is essentially a 
consequence of Stokes's formula). Namely, if 
$W$ is an integral algebraic $q$-dimensional subscheme in 
${\bf A}^n_\C$ (with completion ${\cal W}$ in ${\rm Proj}\, \C[X_0,...,X_n]$) 
and $P_1,...,P_q$ are $q$ polynomials in $n$ variables 
such that $|{\cal W}|\cap \{{}^h P_1=\cdots={}^h P_q=0\}$ is 
finite and included in $\C^n$, then ($[W]$ being understood as the 
integration current free of multiplicities), 
$$
{\rm Res}\left[\matrix{[W]\wedge Q(X) dX_{i_1}\wedge \cdots \wedge dX_{i_q} 
\cr P_1,...,P_q}\right]=0 
$$ 
when ${\rm deg}\, Q \leq \sum_{j=1}^q {\rm deg}\, P_j -q -1$ 
(corollary 3.1) for any ordered subset $\{i_1,...,i_q\}\subset \{1,...,n\}$.  
\vskip 2mm
\noindent
On the other hand, the transformation law holds for our analytic restricted 
residue (see corollary 2.2). Such a transformation law remains valid (in its 
local formulation) for restricted residue symbols defined through the algebraic
 approach 
(see theorem 2.4 in \cite{huku0:gnus}). 
\vskip 2mm
\noindent
Finally, the local residue symbol 
$$
{\rm Res}_{{\cal W};{\cal D}_1,...,{\cal D}_q,\alpha} \, (\omega) 
={\rm Res}_{\C({\cal X})/\C,\alpha} \, \Bigg( \Bigg[ \matrix{\eta \cr 
f_1,...,f_q }\Bigg]\Bigg)\,, 
$$
where 
$$
\omega = {\eta \over f_1 \cdots f_q}\,,
$$
$\eta \in \omega^q_{\C({\cal X})/\C,\alpha}$ and $f_j=0$, $j=1,...,q$, 
is a local equation 
for ${\cal D}_j$ at $\alpha$, equals to $0$ as soon as $\eta=f_j \widetilde \eta$ for 
some $\widetilde \eta \in \omega^q_{\C({\cal X})/\C,\alpha}$
(see also \cite{huku0:gnus}, section 2). The same annihilation property is 
satisfied by the restricted residual current (proposition 2.3). 
\vskip 2mm
\noindent
Our goal in this section is to profit from the fact that both restricted residual objects 
(defined through the algebraic or analytic approach) satisfy the transformation  law,  
the residue formula, the annihilation property, 
in order to show that they coincide. Therefore, we are able to give an 
algebraic formulation of the Proposition 3.1, which is the Theorem 1.1 stated in our 
introduction.   
\vskip 2mm
\noindent
In order to do that, we will need the following technical lemma~: 
\vskip 2mm
\noindent
\begin{lemma} Let $|W|$ be an irreducible $q$-dimensional algebraic set in $\C^n$ and 
$|{\cal W}|$ its Zariski closure in $\P^n(\C)$. Let also 
$P_1,...,P_q$ be $q$ polynomials in 
$\C[X_1,...,X_n]$ such that $V(P)\cap |W|$ is a discrete (hence finite) 
algebraic subset of $\C^n$, with $0 \in V(P)\cap |W|$. Then, there exists $N_0>0$ such 
that, for any integer $N\geq N_0$, one can find $q n+1$ complex parameters 
$u_{jk}$, $ j=1,\dots,q$, $k=1,\dots, n$, $t\in \C^*$, so that, if 
$$
\widetilde P_{j}^{(N,u,t)} (X):=  tP_j(X) +  \Big(\sum\limits_{k=1}^n u_{jk} X_k\Big)^N \,,\quad j=1,...,q\,,
$$
one has~:       
\begin{itemize} 
\item any point $\alpha\in |W|\cap V(\widetilde P^{(N,u,t)})$ but $0$ belongs to $|W|_{\rm reg}$~; 
\item the set 
$$
|{\cal W}|\cap \{[\zeta_0\,:\,...\,:\, \zeta_n]\in \P^n(\C)\,:\, {}^h P^{(N,u,t)}_j (\zeta_0,...,\zeta_n)=0\,,\ j=1,...,q\} 
$$
is contained in $\C^n$.  
\end{itemize} 
\end{lemma} 

\noindent
{\bf Proof.} Since $|W|$ is irreducible and $q$-dimensional, one has ${\rm dim}\, |W|_{\rm sing} <q$~; 
one can 
find an algebraic affine hypersurface $H:=\{\zeta\in \C^n\,;\, H(\zeta)=0\}$ (with Zariski 
closure $|{\cal H}|$) such that $|W|_{\rm sing} \subset H$ 
and ${\rm dim}\, (|{\cal W}|\cap |{\cal H}|)<q$. 
\vskip 1mm
\noindent
Let $N_0 > \deg P_j$, $j=1,...,q$, and $N\geq N_0$. Assume also that 
$N\geq \rho_{P,W} (0)$, where $\rho_{P,W}(0)$ is the order of vanishing of $P$ at the 
origin (along $|W|$). 

\noindent
Let $u=[u_{jk}]$, $j=1,...,q$, $k=1,...,n$, be a $(q,n)$ matrix with generic complex entries,  
$$
M_u:=\{\zeta\in \C^n\,:\, u_{j1}\zeta_1+\cdots + u_{jn} \zeta_n=0\,,\ j=1,...,q\} 
$$
and $|{\cal M}_u|$ its Zariski closure in $\P^n(\C)$. Since ${\rm dim}\, |{\cal W}|=q$ 
and ${\rm dim}\, (|{\cal W}|\cap |{\cal H}|)<q$, 
$|{\cal M}_u| \cap |{\cal W}|\subset \C^n$ and $|{\cal M}_u|\cap |{\cal W}|\cap |{\cal H}|=\{0\}$ 
for $u$ generic. Therefore, 
for such a generic choice of $u$ ($u=u^{0}$) (this choice will be refined 
later), for any $t\in \C^*$, the polynomials 
$$
t P_j(X) + (u^{0}_{j1} X_1+\cdots + u_{jn}^0 X_n)^N \,, \quad j=1,...,q\,,
$$
define in $\C^n$ an algebraic set $Z^{(N,u^0,t)}$ whose closure ${\cal Z}^{(N,u^0,t)}$ in 
$\P^n(\C)$ intersects $|{\cal W}|$ only at points in $\C^n$ (note that $0$ is one 
of these points). The algebraic 
set $|W|\cap Z^{(N,u^0,t)}$ can be described as 
$$
|W|\cap Z^{(N,u^0,t)}=\{\zeta^{(N,1)}(u^0,t),...,\zeta^{(N,m)} (u^0,t)\} \cup \{0\}\,,
$$
where $m$ is fixed (depending on $N$ and $|{\cal W}|$) 
and the $t\mapsto \zeta^{(N,j)} (u^0,t)$, $j=1,...,m$, 
are algebraic $\C^n$-valued functions of $t$ which are not identically $0$ and 
can be classified in two classes, depending on their behavior when 
$|t|$ tends to zero. A branch $t\mapsto \zeta^{(N,j)} (u^0,t)$ will be in the first class if 
$\zeta^{(N,j)} (u^0,t)$ tends to zero when $|t|$ tends to $0$. It will be 
in the second class if $\zeta^{(N,j)} (u^0,t)$ tends to a point in $|W|\cap M_{u^0}$ which is distinct from $0$ 
when $|t|$ goes to $0$. It follows then from $M_{u^0}\cap |W|\cap H=\{0\}$ 
that none of the functions 
$$
t\mapsto H(\zeta^{(N,j)} (u^0,t)) 
$$
where $t\mapsto \zeta^{(N,j)} (u^0,t)$ belongs to the second category, can be identically 
equal to $0$. The behavior of branches of the first category can now be studied 
when $|t|$ goes to infinity. The assumption on $N$ ensures us that such branches either 
approach points in $(|W|\cap V(P)) \setminus \{0\}$, either satisfy
$$
\lim\limits_{|t|\rightarrow \infty} |\zeta^{(N,j)} (u^0,t)|=+\infty 
$$
in the second alternative. The hypothesis on $u^0$ implies that the function 
$t\mapsto H(\zeta^{(N,j)} (u^0,t))$ is not identically $0$ if 
we are in the second alternative. If $u^0$ is conveniently 
choosen (in terms of the Taylor developments at the first order for 
$P_1,...,P_q$ at the points in $(|W|\cap V(P)) \setminus \{0\}$, the 
assertion $t\mapsto H(\zeta^{(N,j)} (u^0,t))\not\equiv 0$ also holds for branches 
concerned by the first alternative. Finally, for any branch $t\mapsto \zeta^{(N,j)} (u^0,t)$, one 
has $H(\zeta^{(N,j)} (u^0,t)\not\equiv 0$. Therefore, once $u^0$ has been 
conveniently chosen, one can pick up $t\not=0$ such that the map  
$\widetilde P^{(N,u^0,t)}$ satisfies the assertions of the lemma. $\quad \diamondsuit$ 

\vskip 2mm
\noindent
We can now relate the analytic and algebraic approaches for restricted residual 
symbols. 

\begin{prop} 
Let ${\cal W}$ be a complete integral $\C$-variety of dimension $q$, embedded in the projective 
scheme ${\rm Proj}\, \C [X_0,...,X_n]$, $\alpha$ be a closed point in $|{\cal W}|$ such 
that $\alpha\in \C^n$ and ${\cal D}_1$,...,${\cal D}_q$ be   
$q$ Cartier divisors on ${\cal W}$ so that the intersection $|{\cal D}_1|\cap \cdots \cap 
|{\cal D}_q|$ define a zero-dimensional scheme on ${\cal W}$ in a neighborhood 
of $\alpha$. If 
$$
\omega = {\eta \over P_1 \cdots P_q}\,,
$$
where $\eta = Q\, dX_{i_1}\wedge \cdots \wedge dX_{i_q}$, 
$Q\in \C[X_1,...,X_n]$ induces an element in $\omega^q_{\C({\cal X})/\C,\alpha}$
 and  
$P_1,...,P_q$ are elements in $\C[X_1,...,X_n]$ such that 
$P_j$, $j=1,...,q$, is a local equation 
for ${\cal D}_j$ at $\alpha$, then, for any function $\varphi \in{\cal D}(\C^n)$ with 
arbitrary small support around $\alpha$ satisfying $\varphi\equiv 1$ in 
a neighborhood of $\alpha$, one has 
\begin{eqnarray}
{\rm Res}_{{\cal W};{\cal D}_1,...,{\cal D}_q,\alpha} \, (\omega) = 
{\rm Res}\, \left[ \matrix{ [W]\wedge \varphi \, \eta \cr 
P_1,..., P_q }\right]\, .  
\end{eqnarray} 

\end{prop} 

\noindent
{\bf Proof.} One can assume for the sake of simplicity that $\alpha=0$. Let ${\cal M}$ 
be the maximal ideal $(X_1,...,X_n)$ in the local algebra ${\cal O}_{\C[X_1,...,X_n],0}$ 
and $(I(W))_0$ the localization at $0$ of the radical ideal 
$$
I(W):=\{g\in \C[X_1,...,X_n]\,;\, 
g(\zeta)=0 \ \forall \zeta \in |{\cal W}|\cap \C^n\}\, .
$$
Choose $p\in \N^*$ such that 
$$
{\cal M}^{p} \in \Big([(P_1,...,P_q)_{0}]^2 , I(W)_0 \Big)\, . 
$$
It follows from the validity of the transformation law and the annihilating 
property in the algebraic context that, if 
$$
\widetilde P_j(X):=P_j(X)+ \Big( \sum\limits_{k=1}^n u_{jk} X_k \Big)^{p},
$$
 then one has, for any $\eta= Q\, dX_{i_1}\wedge \cdots \wedge dX_{i_q}$, 
$Q\in \C[X_1,...,X_n]$, that
\begin{eqnarray}
{\rm Res}_{{\cal W};{\cal D}_1,...,{\cal D}_q,0} \, \Big( {\eta \over P_1 \cdots P_q}\Big) 
={\rm Res}_{{\cal W};\widetilde {\cal D}_1,...,\widetilde {\cal D}_q,0} \, \Big( {\eta \over \widetilde P_1 \cdots 
\widetilde P_q}\Big)\,,
\end{eqnarray}
where $\widetilde {\cal D}_j$, $j=1,...,q$, is the Cartier divisor on ${\cal W}$ with local equation 
$\widetilde P_j$ in a neighborhood of the origin. On the 
other hand, it follows from Proposition 2.2 and Proposition 2.3 
that, for any test-function $\varphi$ with arbitrary small support around the origin, one has 
also
\begin{eqnarray}
{\rm Res}\, \left[ \matrix{ [W]\wedge \varphi \, \eta \cr 
P_1,..., P_q }\right] = 
{\rm Res}\, \left[ \matrix{ [W]\wedge \varphi\, \eta \cr 
\widetilde P_1,..., \widetilde P_q }\right]\, . 
\end{eqnarray}
If the $u_{jk}$, $j=1,...,q$, $k=1,...,n$ are generic
(see for example the construction in the proof of lemma 4.1), the algebraic 
set $V(\widetilde P)\cap |{\cal W}|\cap \C^n$ is discrete (hence finite).  
We can then conclude from (4.2) and (4.3) that in order to prove (4.1), it is 
not restrictive to assume that the algebraic set $V(P)\cap |{\cal W}| \cap \C^n$ 
is finite, what we will do from now on.   
\vskip 2mm
\noindent
The same argument as above shows that, in order to prove (4.1), one can replace 
$P_j$, $j=1,...,q$, by the polynomial  
$$
{1\over t} \widetilde P_j^{(N,u,t)}   
$$
constructed in lemma 4.1 ($N$ being choosen sufficiently large, certainly such that 
$N\geq \max \deg P_j$, ${\cal M}^N \subset (I(P)_0,I(W)_0)$) 
and $\deg Q < q (N-1)$), and this is what we do 
(preserving the notations $P_j$ and ${\cal D}_j$).  
As a consequence 
of the residue formula in the algebraic context (which we recalled at the beginning of this 
section) and of Corollary 3.1, one has 
\begin{eqnarray}
&&\sum\limits_{\alpha \in V(P) \cap {\cal W} (\C)} 
{\rm Res}_{{\cal W};{\cal D}_1,...,{\cal D}_q,\alpha} \, \Big( {\eta \over P_1 \cdots P_q}\Big) 
=\sum\limits_{\alpha \in V(P) \cap {\cal W} (\C)} 
 {\rm Res}\, \left[ \matrix{ [W]\wedge \varphi \, \eta \cr 
P_1,..., P_q }\right]\nonumber \\
&&\end{eqnarray}
whenever $\varphi$ is a test-function in ${\cal D}(\C^n)$ with
 arbitrary small  support
around the points $\alpha \in V(P)\cap |{\cal W}|$, such that $\varphi\equiv 1$ in 
a neighborhood of each of these points ($\varphi_\alpha$ will 
denote next $\varphi\,\theta_\alpha$, where $\theta_\alpha$ is a test-function
with support arbitrary small around $\alpha$ and $\theta_\alpha\equiv 1$ in a 
neighborhood of $\alpha$). If $\alpha$ is any point in $V(P)\cap |{\cal W}|$ 
distinct from $0$, ${\cal W}$ is smooth about $\alpha$ 
(lemma 4.1, first assertion) and we know in this case 
that  
\begin{eqnarray*}
{\rm Res}_{{\cal W};{\cal D}_1,...,{\cal D}_q,\alpha} \, \Big( {\eta \over P_1 \cdots P_q}\Big) 
&=&{\rm Res}\, \left[ \matrix{ [W]\wedge \varphi_\alpha\,  \eta \cr 
P_1,..., P_q }\right]\,,  
 \end{eqnarray*}
since the construction of our restricted residual currents corresponds to 
the construction proposed in \cite{gh:gnus}, chap 5, section 1 (this is a consequence of the 
classical relation between Bochner-Martinelli and Cauchy kernels), which is known to 
fit with the algebraic approach in the smooth case (as it was recalled at the beginning of this section).
Formula (4.1) follows then from (4.4) and from the identifications (4.5). $\quad\diamondsuit$  

\vskip 2mm
\noindent
{\bf Proof of Theorem 1.1.} 
We may now transpose to the algebraic context the analytic result stated in proposition 3.1. 
This gives the statement of the Theorem 1.1 of our introduction, provided we remember that 
we have 
$$
{\rm Res}_{{\cal W};{\cal D}_1,...,{\cal D}_q,\alpha} \, (\omega) 
={\rm Res}_{\C({\cal W})/\C,\alpha} \, \Bigg( \Bigg[ \matrix{\eta \cr 
f_1,...,f_q }\Bigg]\Bigg) 
$$
for any point $\alpha$ in $|{\cal W}|\cap |{\cal D}_1|\cap 
\cdots \cap |{\cal D}_q| \cap \C^n$ (here we just assume 
that ${\cal D}_1$,...,${\cal D}_n$ define a $0$-dimensional 
scheme on $W$, there is no assumption about what happens on 
$|{\cal W}|\setminus |W|$)  
and any $\omega$ in $\Omega^q_{\C({\cal W})/\C}$ 
with simple poles (in $W$) along ${\cal D}_1+\cdots +{\cal D}_q$ ($\eta=f_1\cdots f_q \,\omega$, 
where $f_j$ denotes a local equation for the Cartier divisor ${\cal D}_j$).
Since the reference to the divisors ${\cal D}_1,...,{\cal D}_q$ was implicit in 
the expression of the element in $\Omega^q_{\C({\cal W})/\C}$, 
we used the abridged notation ${\rm Res}_{W,\alpha} [\ ]$ instead     
of ${\rm Res}_{{\cal W};{\cal D}_1,...,{\cal D}_q,\alpha}$ 
in order to formulate the statement in this theorem. $\diamond $\\
As a direct consequence we formulate the restricted version of the 
Cayley-Bacharach Theorem.
\begin{corollary} 
Let $W$ be a $q$-dimensional irreducible affine algebraic subvariety 
in ${\bf A}^n_\C$ \($0<q<n$\) 
and $P_1,...,P_q$ be $q$ polynomials in $\C[X_1,...,X_n]$ satisfying
 the condition
$(1.2)$. Assume also that $V(P)$ and $\vert W\vert $ intersect transversally
at any of the $k$ points in $V(P)\cap \vert W\vert $. 
Then any 
algebraic hypersurface $\{Q=0\}$, $Q\in \C[X_1,...,X_n]$, 
such that ${\rm deg}\, Q <\delta_1+\cdots +\delta_q -q$,
 which passes through any
$k-1 $ points of the set $V(P)\cap W $ passes through the last one also.
\end{corollary}  

\section{An affine version of Wood's theorem.} 

Let $\gamma_1$, ...,$\gamma_d$ be $d$ pieces of manifold in 
$\P^n(\C)$ and $|{\cal L}_{0,0}|$ be a line in $\P^{n}(\C)$ which 
intersects each of the $\gamma_j$ transversally respectively 
at distinct points $p_{j0}$, $j=1,...,d$. Assume that affine coordinates are 
such that the support $|{\cal L}_{0,0}|$ is the line $\zeta_1=\cdots=\zeta_{n-1}=0$. 
Then, for $(\alpha,\beta)\in (\C^{n-1})^2$ close to $(0,0)$, the 
projective line 
$$
|{\cal L}_{\alpha,\beta}|:=\{[\zeta_0:\,... \, :\zeta_n]\in \P^n(\C)\,;\, 
\zeta_k= \alpha_k\, \zeta_n +\beta_k\, \zeta_0\,,\ 
k=1,...,n-1\} 
$$
intersects tranversally $\gamma_1$,...,$\gamma_d$ at 
the respective points $p_1(\alpha,\beta)$,..., $p_d(\alpha,\beta)$ ($p_j(\alpha,\beta)$ 
being close to 
$p_{j0}$). In 
\cite{wood1:gnus}, J. Wood gave a simple criterion for the local germs of manifold   
$\gamma_1,...,\gamma_d$ to be germs of a global algebraic hypersurface 
(with degree $d$) $|{\cal H}|$ in $\P^n(\C)$ satisfying the relation 
such that 
$$
|{\cal H}|\cap |{\cal L}_{\alpha,\beta}| =
\{p_1(\alpha,\beta),...,p_d(\alpha,\beta)\}
$$ 
for $(\alpha,\beta)$ close to $(0,0)$. 
The (necessary and sufficient) condition he gave can be formulated as 
follows~: 
\begin{eqnarray}
\sum\limits_{j=1}^d \zeta_n [p_j(\alpha,\beta)]= h_0(\alpha)+\sum\limits_{k=1}^{n-1} 
h_k(\alpha)\, \beta_k\,,
\end{eqnarray}
where $h_0,...,h_{n-1}$ are germs of holomorphic functions in $\alpha$ 
at the origin (here $\zeta_n[p]$, where $p$ denotes a point in $\C^n$,  
means the $n$-th affine coordinate of $p$). Note that the 
algebraic hypersurface $|{\cal H}|$ (in $\P^n(\C)$) which interpolates $\gamma_1,...,\gamma_d$ 
is such that 
its intersection at infinity with any line $|{\cal L}_{\alpha,\beta}|$, with $(\alpha,\beta)$ 
close to $(0,0)$, is empty. 
What we would like to state here is an affine analog of this result, ${\bf P}^n(\C)$ being replaced 
by some irreducible $q$-dimensional affine algebraic subvariety of $\C^n$ ($q=2,...,n$).       

\vskip 2mm
\noindent
Let us first state the following easy consequence of our Theorem 1.1.
\vskip 2mm 
\begin{prop} Let $W$ be an algebraic irreducible $q$-dimensional subvariety of the 
affine scheme ${\bf A}^n_\C$ (with $2\leq q\leq n$), 
$m$ be a positive integer strictly between $0$ and $q$,  
and $\gamma_1$,...,$\gamma_d$ be $d$ disjoint pieces of $q-m$-dimensional 
analytic manifold such that $\gamma_j$ lies in $|W|_{\rm reg}$ for $j=1,...,d$.
 Furthermore, assume that the affine 
$n+m-q$-dimensional subspace 
$$
L_{0,0}:=\{\zeta \in \C^n\,;\, \zeta_k=0\,,\ k=1,...,q-m\}
$$
intersects each $\gamma_j$ transversally respectively 
at points $p_{j0}$, $j=1,...,d$. Suppose that there are strictly 
positive rational numbers $\delta_1,...,\delta_m$ and polynomials $P_1,...,P_m$ with 
$\deg P_j=d_j\geq \delta_j$, $j=1,...,m$, such that 
\begin{itemize} 
\item $|W| \cap V(P)$ is a $q-m$-dimensional variety in $\C^n$ which interpolates 
the pieces $\gamma_j$ and is such that $|W|\cap V(P)\cap L_{0,0}=\{p_{10},...,p_{d0}\}$~; 
\item there exists strictly positive constants $\kappa, K$ such that 
\begin{eqnarray}
\zeta \in |W|\,,\ \|\zeta\|\geq K \Longrightarrow
\sum\limits_{j=1}^m {|P_j(\zeta)|\over \|\zeta\|^{\delta_j}} + 
\sum\limits_{k=1}^{q-m} {|\zeta_k|\over \|\zeta\|} \geq \kappa\, . 
\end{eqnarray} 
\end{itemize}  
Then, for $(\alpha,\beta)$ close to $(0,0)$ in $(\C^{n+m-q})^{q-m} \times \C^{q-m}$, 
the affine $n+m-q$-dimensional subspace 
$$
L_{\alpha,\beta}:=\Big\{\zeta \in \C^n\,;\, \zeta_k=\sum\limits_{r=1}^{n+m-q} 
\, \alpha_{k,r} \zeta_{q-m+r} + \beta_k\,,  
\ k=1,...,q-m\Big\}
$$
intersects each $\gamma_j$ transversally respectively at the points 
$p_j(\alpha,\beta)$, $j=1,...,d$ (necesseraly distinct and  
close to the $p_{j0}$) and one has 
\begin{eqnarray}
\sum\limits_{j=1}^d \zeta_l [p_j(\alpha,\beta)]= 
\sum\limits_{{\underline k\in \N^{q-m}}\atop {|\underline k|
\leq \rho +1}} 
h_{\underline k}^{(l)} (\alpha)\,  \beta_1^{k_1}\cdots \beta_{q-m}^{k_{q-m}}\,,\ 
l=q-m+1,...,n\,,
\end{eqnarray}
where the $h_{\underline k}^{(l)}$ are germs of holomophic functions in $\alpha$ about the origin 
and 
$$
\rho:=\sum_{j=1}^m (d_j-\delta_j)
$$  
\end{prop}  

\noindent
{\bf Proof.} Let, for $k=1,...,q-m$, 
$$
\Lambda_{\alpha,\beta,k}(\zeta):=\zeta_k-\sum\limits_{r=1}^{n+m-q} \alpha_{k,r} \, 
\zeta_{q-m+r} -\beta_k\,,\ \zeta \in \C^n \, ;
$$   
condition $(5.7)$ implies that, when $(\alpha,\beta)$ is 
sufficiently close to $(0,0)$, one has 
\begin{eqnarray}
\zeta \in |W|\,,\ \|\zeta\|\geq K \Longrightarrow
\sum\limits_{j=1}^m {|P_j(\zeta)|\over \|\zeta\|^{\delta_j}} + 
\sum\limits_{k=1}^{q-m} {|\Lambda_{\alpha,\beta,k} (\zeta)|\over \|\zeta\|} 
\geq {\kappa\over 2}\, . 
\end{eqnarray} 
This shows that for $(\alpha,\beta)$ close to $(0,0)$, the only points 
in $L_{\alpha,\beta} \cap |W| \cap V(P)$ are $d$ points $p_j(\alpha,\beta)$, 
$j=1,...,d$ which approach the points $p_{10},...,p_{d0}$ (about each of 
these points, one can use the implicit function theorem in order to 
describe the intersection $\gamma_j \cap L_{\alpha,\beta}$). This proves the 
first assertion of the proposition.
\vskip 1mm
\noindent  
It follows from proposition 3.1 that, as soon as 
the multi-index $\underline k\in \N^{q-m}$ is 
such that 
$$
\sum\limits_{j=1}^m (d_j-1) +1 < \sum\limits_{j=1}^m \delta_j + 
\sum\limits_{l=1}^{q-m} (k_l+1)-q=\sum\limits_{j=1}^m \delta_j+ 
|\underline k|-m\,,
$$
then, for $l=q-m+1,...,n$, for any finite ordered subset $\{i_1,...,i_{q-m}\}
\subset \{1,...,n\}$,   
$$
{\rm Res}\, \left[\matrix{ [W] \wedge X_l \, \Big(\bigwedge\limits_{j=1}^m 
dP_j\Big) \wedge \Big(\bigwedge\limits_{l=1}^{q-m} dX_{i_l}\Big) \cr \cr 
P_1,...,P_m, (\Lambda_{\alpha,\beta,1})^{k_1+1},..., 
(\Lambda_{\alpha,\beta,q-m})^{k_{q-m}+1} }\right]=0 
$$
for $(\alpha,\beta)$ such that (5.9) holds. It is immediate to check 
(use for example formula $(3.7)$) that 
for such $(\alpha,\beta)$, one has, for any multi-index 
$\underline k=(k_1,...,k_{q-m})\in \N^{q-m}$, 
\begin{eqnarray}
&& {\partial^{|\underline k|} \over 
\partial \beta_1^{k_1}\cdots 
\partial \beta_{q-m}^{k_{q-m}}} 
\, {\rm Res}\, \left[\matrix{ [W] \wedge X_l \, \Big(\bigwedge\limits_{j=1}^m 
dP_j\Big) \wedge \Big(\bigwedge\limits_{l=1}^{q-m} dX_{i_l}\Big) \cr \cr 
P_1,...,P_m, \Lambda_{\alpha,\beta,1},..., 
\Lambda_{\alpha,\beta,q-m}}\right] \nonumber \\
&& 
\quad =\pm\,  {\rm Res}\, \left[\matrix{ [W] \wedge X_l \, \Big(\bigwedge\limits_{j=1}^m 
dP_j\Big) \wedge \Big(\bigwedge\limits_{l=1}^{q-m} dX_{i_l}\Big) \cr \cr 
P_1,...,P_m, (\Lambda_{\alpha,\beta,1})^{k_1+1},..., 
(\Lambda_{\alpha,\beta,q-m})^{k_{q-m}+1} }\right]\, . 
\end{eqnarray} 
Then it follows from $(5.9)$ that the right-hand side of $(5.10)$ 
(hence the left-hand side) equals 
identically $0$ when 
$$
|\underline k| > \sum\limits_{j=1}^m (d_j-\delta_j) +1=
\rho+1\,.  
$$
This proves that, when $(\alpha,\beta)$ is close to 
$(0,0)$ and $l=m-q+1,...,n$,    
$$
\sum\limits_{j=1}^d \zeta_l [p_j(\alpha,\beta)] 
\equiv {\rm Res}\, \left[\matrix{ [W] \wedge \zeta_l \, \Big(\bigwedge\limits_{j=1}^m 
dP_j\Big) \wedge \Big(\bigwedge\limits_{l=1}^{q-m} d\Lambda_{\alpha,\beta,l}
\Big) \cr \cr 
P_1,...,P_m, \Lambda_{\alpha,\beta,1},..., 
\Lambda_{\alpha,\beta,q-m}}\right]
$$
is a polynomial expression in $\beta=(\beta_1,...,\beta_{q-m})$ with 
total degree at most $\rho+1$ (the coefficients being holomorphic 
functions in $\alpha$). The second assertion of the proposition is proved. 
$\quad\diamondsuit$ 
\vskip 2mm
\noindent
{\bf Remark.} Note that we recover here as a particular case 
the necessity of Wood's condition in the 
case $W={\bf A}^n_\C$, $m=1$, $\delta_1=d_1=d$, 
which means precisely that in this case we also impose 
the restriction 
$$
\{\widetilde \zeta \in \P^n(\C)\,;\, {}^h P_1(\widetilde \zeta)=0\}\cap |{\cal L}_{0,0}|=
\{p_{10},...,p_{d0}\}\,.
$$ 
\vskip 2mm
\noindent
Furthermore, one can state the following proposition, which appears as 
a weak converse of proposition 5.1 in the affine setting. 
\vskip 2mm
\noindent
\begin{prop} 
Let $\gamma_1$,...,$\gamma_d$ be $d$ disjoint pieces of $n-m$-dimensional 
analytic manifold ($1\leq m <n$) in the affine space $\C^n$.  
Suppose that for any $(\alpha,\beta)\in (\C^{m})^{n-m} \times \C^{n-m}$ 
close to $(0,0)$, the affine $m$-dimensional subspace 
$$
L_{\alpha,\beta}:=\Big\{\zeta \in \C^n\,;\, \zeta_k=\sum\limits_{r=1}^{m} 
\, \alpha_{k,r} \zeta_{n-m+r} + \beta_k\,,  
\ k=1,...,n-m\Big\} 
$$
intersects transversally $\gamma_1$,...,$\gamma_d$ respectively at points 
$p_1(\alpha,\beta),...,p_d(\alpha,\beta)$. 
Assume that there exists $D\in \N$ and analytic 
functions $h_{\underline k}^{(l)}$, $|\underline k|\leq D+1$, $l=n-m+1,...,n$,  
in a neighborhood of $0$ in $(\C^{m})^{n-m}$ such that for 
$(\alpha,\beta)$ close to $(0,0)$ in $(\C^{m})^{n-m} \times \C^{n-m}$, 
for any $l=n-m+1,...,n$, 
\begin{eqnarray}
\sum\limits_{j=1}^d \zeta_l [p_j(\alpha,\beta)]= 
\sum\limits_{{\underline k\in \N^{n-m}}\atop {|\underline k|
\leq D+1}} 
h_{\underline k}^{(l)} (\alpha)\,  \beta_1^{k_1}\cdots \beta_{q-m}^{k_{n-m}}\, .
\end{eqnarray}
Then, one can find a collection of 
polynomials $(P_\iota)_{\iota \in {\cal J}}$ 
with degree at most $d+D$ which define an affine algebraic variety $V(P)$ 
such that for some convenient constants $\epsilon>0$, $\kappa>0$, $K>0$, one 
has~:  
\begin{itemize} 
\item if 
$$
\Gamma_\epsilon:= 
\bigcup\limits_{{(\alpha,\beta) \in (\C^{m})^{n-m}\times \C^{n-m}}\atop 
{\max (\|\alpha\|,\|\beta\|) <\epsilon}} L_{\alpha,\beta}\,,
$$
then
\begin{eqnarray}
\zeta \in \Gamma_\epsilon\,,\ \|\zeta\| \geq K \ \Longrightarrow 
\max\limits_{\iota \in {\cal J}} |P_\iota (\zeta)| \geq \kappa \|\zeta\|^d\,; 
\end{eqnarray} 
\item for $\max (\|\alpha\|,\|\beta\|) <\epsilon$, one has 
\begin{eqnarray}
L_{\alpha,\beta} \cap V(P)=\{p_1 (\alpha,\beta),...,p_d(\alpha,\beta)\}\, .
\end{eqnarray}
\end{itemize} 
\end{prop}
Before the proof of this 
proposition we remark that our approach is directly inspired from 
\cite{wood1:gnus} (page 237, 
proof of the sufficiency). First we observe, as in Wood's argument,  
that conditions (5.11) imply 
that for any integer $\sigma\in \N^*$, for any $l=n-m+1,...,n$, one has, for 
$(\alpha,\beta)$ close to $(0,0)$ in $(\C^{m})^{n-m} \times \C^{n-m}$,  
\begin{eqnarray}
\sum\limits_{j=1}^d (\zeta_l[p_j(\alpha,\beta)])^\sigma  
=\sum\limits_{{\underline k\in \N^{n-m}}\atop {|\underline k|
\leq D+\sigma}} 
h_{\sigma,\underline k}^{(l)} (\alpha)\,  \beta_1^{k_1}\cdots \beta_{n-m}^{k_{n-m}}\,,
\end{eqnarray}
where the $h_{\sigma,\underline k}^{(l)}$ are analytic functions in $\alpha$ in a neighborhood 
of $0$.\\
 Actually, if $\gamma_j$ is defined semi-locally as the 
smooth complete intersection 
$$
\gamma_j= \{\zeta\in \C^n\,:\, \Phi_{j,1}(\zeta)=\cdots=\Phi_{j,m}(\zeta)=0\}\,,
$$
then  for any $j=1,...,d$ and  for any $k=1,...,n-m$,     
\begin{eqnarray*} 
\zeta_k (p_j(\alpha,\beta))&=&{\rm Res} \left[\matrix{\varphi_j \, \zeta_k\,  
d\Phi_j \wedge d\Lambda_{\alpha,\beta} 
\cr \Phi_{j,1},...,\Phi_{j,m},\Lambda_{\alpha,\beta,1},..., \Lambda_{\alpha,\beta, n-m}}
\right]\\   
&=&   
{\rm Res} \left[\matrix{\varphi_j \, \Big(\sum\limits_{r=1}^{m} \alpha_{k,r} 
\zeta_{n-m+r}+\beta_k\Big)\, 
d\Phi_j \wedge d\Lambda_{\alpha,\beta}
\cr \Phi_{j,1},...,\Phi_{j,m},\Lambda_{\alpha,\beta,1},..., \Lambda_{\alpha,\beta, n-m}}
\right]\,,
\end{eqnarray*}
where 
$$
\Lambda_{\alpha,\beta,k}:=\zeta_k-\sum\limits_{r=1}^m 
\alpha_{k,r}\, \zeta_{n-m+r} -\beta_k\,,\quad k=1,...,n-m\,, 
$$ 
$$
d\Phi_j:= \bigwedge_{k=1}^{m} d\Phi_{j,k}\quad ,\quad d\Lambda_{\alpha,\beta}:= 
\bigwedge_{k=1}^{n-m} d\Lambda_{\alpha,\beta,k}\,,
$$
and $\varphi_j$ is a test-function which 
is supported by an arbitrary small neighborhood of $p_{j0}$ and is such that 
$\varphi\equiv 1$ near this point. 
This implies that one can also represent in fact 
any sum $\sum_j \zeta_l[p_j(\alpha,\beta)]$\,,
$l=1,...,n$,  
in the form (5.11). Since we also have for any $j=1,...,d$, 
for any $\sigma \in \N^*$ and for any $k=1,...,n-m$ that 
$$
{\rm Res} \left[\matrix{\varphi_j \, \Big(\zeta_k-\sum\limits_{r=1}^{m} \alpha_{k,r} 
\zeta_{n-m+r}-\beta_k \Big)^{\sigma}\,  
d\Phi_j \wedge d\Lambda_{\alpha,\beta}
\cr \Phi_{j,1},...,\Phi_{j,m},\Lambda_{\alpha,\beta,1},..., \Lambda_{\alpha,\beta, n-m}}
\right] =0\,,  
$$
it follows by induction on $\sigma$ that (5.14) holds in a 
neighborhood of $(0,0)$ for any $l=1,...,n$, in particular for any 
$l=n-m+1,...,n$ 
(here we use the fact that one can choose $\alpha$ generic in a neighborhood of 
the origin in $(\C^m)^{n-m}$). At this point we return to the proof of the 
Proposition 5.2.
\\
\noindent
\noindent
{\bf Proof of Proposition 5.2.}   
For any $l=n-m+1,...,n$,  
let $A_l$ be the polynomial in variable $X_l$ (with coefficients analytic in $\alpha$ and 
polynomial in $\beta$) defined as 
\begin{eqnarray}
A_l(X_l,\alpha\,;\, \beta) &=& \prod\limits_{j=1}^d (X_l-\zeta_l[p_j(\alpha,\beta)]) \nonumber \\ 
&=& X_l^d - A_{l1}(\alpha,\beta) X_l^{d-1} + \cdots + (-1)^d A_{ld}(\alpha,\beta) \nonumber \\
\end{eqnarray}
($\alpha$ and $\beta$ close to $0$ in their respective spaces). For any such $\alpha$, 
denote as $P_{l,\alpha}$ the element in $\C[X_1,...,X_n]$ defined as 
$$
P_{l,\alpha} (X)=A_l\Big(X_l,\alpha\,;\, X_{1}-\sum\limits_{
r=1}^{m} \alpha_{1,r} X_{n-m+r}, ..., X_{n-m}-\sum\limits_{r=1}^{m} 
\alpha_{n-m,r} X_{n-m+r}\Big)\, .  
$$
For each $\alpha$ close to $0$, $P_{l,\alpha}$ is a polynomial in variables
$(X_1,...,X_n)$ with total degree less than $d+D$. If $\zeta$ is 
a point in $\gamma_j$, then one can write, for any $\alpha$ close to 
$0$ in $(\C^{m})^{n-m}$, for any $l=n-m+1,...,n$, 
$$
\zeta_l= \zeta_l [p_j(\alpha,\beta_\zeta)]\,,
$$
where 
$$
\beta_{\zeta,k}=\zeta_k -\sum\limits_{r=1}^{m} \alpha_{k,r}\, \zeta_{n-m+r}\,;
$$
then, it follows from the definition (5.15) of $A_l$, hence of $P_{l,\alpha}$, that 
$\zeta \in V(P_{l,\alpha})$ for any $l=n-m+1,...,n$ and for any $\alpha$ 
close to zero. All pieces of manifold $\gamma_1,...,\gamma_d$ lie in 
$V(P_{n-m+1,\alpha},...,P_{n,\alpha})$ for any $\alpha$ close to 
$0$ in $(\C^{m})^{n-m}$. 
Moreover, if $\epsilon$ is sufficiently small and $K=K_\epsilon>0$ large enough, 
then for any $\alpha$ such that 
$\|\alpha\|<\epsilon$, define the strip 
$$
\Gamma_{\alpha,\epsilon}:= \{\zeta \in \C\,:\, \|\zeta\|\geq K\,,\,   
\max\limits_{k=1,...,n-m} 
\Big| \zeta_k -\sum\limits_{r=1}^{m} \alpha_{k,r}\, \zeta_{n-m+r} \Big|
<\epsilon\}\; . 
$$
For  any $\zeta \in \Gamma_{\alpha,\epsilon} $ one has that 
\begin{eqnarray}
\max\limits_{l=n-m+1,...,n} |P_{l,\alpha} (\zeta)|=
\max\limits_{l=n-m+1,...,n} 
\prod\limits_{j=1}^d \Big|\zeta_l- 
\zeta_l[p_j(\alpha, \beta_\zeta)]\Big| \geq  \kappa \|\zeta\|^d  \,, \nonumber \\
\end{eqnarray}
where $\kappa=\kappa_\epsilon$ is a strictly positive constant (independent of $\alpha$).
\vskip 1mm
\noindent
Let now ${\cal F}$ be the finite subset in $L_{0,0}$ defined as 
$$
\zeta \in {\cal F} \Longleftrightarrow 
\forall \, l=n-m+1,...,n\,,\ 
\exists j\in \{1,...,d\}\,,\ \zeta_l=\zeta_l[p_{j0}]
$$
and ${\cal F}':={\cal F}\setminus \{p_{10},...,p_{d0}\}$. There exists an affine  
form $\Lambda$ in $X_{n-m+1},...,X_n$ such that for any $\zeta \in {\cal F}'$, 
$$
\Lambda(\zeta) \not= \Lambda (p_{j0}), \ j=1,...,d \,.
$$
We can define 
$$
B(X_{n-m+1},...,X_n\,;\, \beta) = \prod\limits_{j=1}^d 
(\Lambda(X_{n-m+1},...,X_n))-\Lambda[p_j(0,\beta)])
$$ 
and $Q(X)=B(X_{n-m+1},...,X_n\,;\, X_1,...,X_{n-m})$. One can check as before that 
$V(Q)$ contains the pieces $\gamma_1,...,\gamma_q$. Moreover $Q$ does not vanish at 
any point in ${\cal F}'$, and $\deg Q \leq d+D$. 
\vskip 1mm
\noindent
We now define the family $(P_{\iota})_{\iota \in {\cal J}}$ as the collection 
of all polynomials $P_{l,\alpha}$, $\|\alpha\|<\epsilon$, $l=n-m+1,...,n$, and the 
polynomial $Q$. It is clear that $(5.12)$ holds (since it already holds for the 
collection $(P_{l,\alpha})$, $\|\alpha \|<\epsilon$, $l=n-m+1,...,n$). Since 
the pieces $\gamma_1,...,\gamma_d$ lie in $V(P)$, one has, for any $(\alpha,\beta)$ 
such that $\max (\|\alpha\|, \|\beta\|)<\epsilon$, 
$$
\{p_1(\alpha,\beta),...,p_d(\alpha,\beta) \} 
\subset L_{\alpha,\beta} \cap V(P)\, .
$$
Since points in $V(P) \cap \Gamma_\epsilon$ lie in the compact 
set $\{\zeta\in \C^n\,:\, \|\zeta\| \leq K\}$ (because of $(5.12)$), points 
in $L_{\alpha,\beta} \cap V(P)$ uniformly approach points in 
$L_{0,0}\cap V(P)$ when $(\alpha,\beta)$ tends to $(0,0)$. Since 
$Q\not=0$ on ${\cal F}'$ and 
$$
V(P_{n-m+1,0},...,P_{n,0})\cap L_{0,0}
\subset {\cal F}\,,
$$
one has that $V(P)\cap L_{0,0}=\{p_{10},...,p_{d0}\}$. Therefore, 
if we refine the choice of $\epsilon$, we can assume that (5.13) 
holds. This concludes the proof of our proposition. $\quad\diamondsuit$    
\vskip 2mm
\noindent
In the particular case $m=1$, one can be more precise and 
repeat Wood's argument in order to obtain the~: 
\vskip 2mm
\noindent
\begin{prop} 
Let $\gamma_1$,...,$\gamma_d$ be $d$ disjoint pieces of smooth  
analytic hypersurface in the affine space $\C^n$.  
Suppose that for any $(\alpha,\beta)\in (\C^{n-1})^2$ 
close to $(0,0)$, the affine line  
$$
L_{\alpha,\beta}:=\Big\{\zeta \in \C^n\,;\, \zeta_k=\alpha_k \zeta_n + \beta_k\,,  
\ k=1,...,n-1\Big\} 
$$
intersects transversally $\gamma_1$,...,$\gamma_d$ respectively at points 
$p_1(\alpha,\beta),...,p_d(\alpha,\beta)$. 
Assume that there exists $D\in \N$ and analytic 
functions $h_{\underline k}$, $|\underline k|\leq D+1$,  
in a neighborhood of $0$ in $\C^{n-1}$ such that for 
$(\alpha,\beta)$ close to $(0,0)$ in $(\C^{n-1})^2$, 
\begin{eqnarray*}
\sum\limits_{j=1}^d \zeta_n [p_j(\alpha,\beta)]= 
\sum\limits_{{\underline k\in \N^{n-m}}\atop {|\underline k|
\leq D+1}} 
h_{\underline k} (\alpha)\,  \beta_1^{k_1}\cdots \beta_{q-m}^{k_{n-m}}\,. 
\end{eqnarray*}
(one from the $h_{\underline k}$ for $|\underline k|=D+1$ being non identically zero).   
Then, one can find a polynomial $P$ 
with degree $d+D$ which defines an affine algebraic variety $V(P)$ 
such that for some convenient constants $\epsilon>0$, $\kappa>0$, $K>0$, one 
has~:  
\begin{itemize} 
\item if 
$$
\Gamma_\epsilon:= 
\bigcup\limits_{{(\alpha,\beta) \in (\C^{m})^{n-m}\times \C^{n-m}}\atop 
{\max (\|\alpha\|,\|\beta\|) <\epsilon}} L_{\alpha,\beta}\,,
$$
then
\begin{eqnarray*}
\zeta \in \Gamma_\epsilon\,,\ \|\zeta\| \geq K \ \Longrightarrow 
|P (\zeta)| \geq \kappa \|\zeta\|^d\,; 
\end{eqnarray*} 
\item for $\max (\|\alpha\|,\|\beta\|) <\epsilon$, one has 
\begin{eqnarray*}
L_{\alpha,\beta} \cap V(P)=\{p_1 (\alpha,\beta),...,p_d(\alpha,\beta)\}\,.
\end{eqnarray*}
\end{itemize} 
\end{prop}

\noindent
{\bf Proof.} We repeat here Wood's argument (as we did in the 
proof of proposition 5.2). For $(\alpha,\beta)$ close to 
$(0,0)$, the points $p_1(\alpha,\beta),...,p_d(\alpha,\beta)$ are 
the only intersection points (in the affine space $\C^n$) between the 
affine line $L_{\alpha,\beta}$ and the affine algebraic hypersurface (with exact degree 
$d+D$ for $\alpha$ generic) $\{P_{n,\alpha}=0\}$. 
Since these intersection points are simple 
(the line hits each $\gamma_j$ transversally), the homogeneous 
polynomial ${}^ h P_{n,\alpha}$ vanishes at the order at most  
$D$ at the point $p_{\alpha,\infty}$ at infinity 
on the projective line $|{\cal L}_{\alpha,0}|$. 
On the other hand, one has 
$$
|P_{n,\alpha} (\zeta)| \geq \kappa \|\zeta\|^d 
$$
in a tube $\Gamma_{\alpha,\epsilon}$ along the line $L_{\alpha,0}$ 
(see (5.16) in the proof of proposition 5.2). This implies that the 
the homogeneous 
polynomial ${}^h P_{n,\alpha}$ vanishes at the order at least  
$D$ at the point $p_{\alpha,\infty}$. Therefore, the 
hypersurface $\{P_{n,\alpha}=0\}$ (with degree $d+D$) contains 
the germs $\gamma_1,...,\gamma_d$ (as simple germs) and the germ corresponding to the 
hyperplane ${\cal H}_\infty$ (this germ being counted with multiplicity $D$). 
This, combined with the fact that the degree of $P_{n,\alpha}$ is exactly $d+D$, implies that 
all the polynomials $P_{n,\alpha}$ define the same algebraic hypersurface 
$H$ in $\C^n$. Since they have the same degree, they are equal (up to some 
constant) to a polynomial $P$. As the auxiliary construction of the polynomial 
$Q$ is not needed in the hypersurface case, proposition 5.3 follows immediately 
from proposition 5.2. $\quad\diamondsuit$ 
\vskip 2mm
\noindent
{\bf Remark.} In the particular case $W={\bf A}^n_\C$, $m=1$, 
proposition 5.3 appears as the 
reciprocal assertion to proposition 5.1. The difficulty in the 
more general case  
$W={\bf A}^n_\C$, $m>1$, is 
to be able to interpolate the germs $\gamma_1,...,\gamma_d$ by an algebraic complete 
intersection $V(P_1,...,P_m)$. It does not seem possible when 
$m>1$ even if 
conditions (5.11) are satisfied with $D=0$ (which would mean that the projective 
variety $\{{}^h P_1=...={}^h P_m=0\}$ 
corresponding to the complete intersection $V(P)$ that 
interpolates the pieces $\gamma_j$ does not hit $|{\cal H}_\infty| \cap |{\cal L}_{0,0}|$).
We do not have the answer to that question yet. Nevertheless, 
proposition 5.2 can be seen as an attempt to settle a converse to proposition 5.1 in general.

\end{document}